\documentclass[preprint,review,3p,numbers,11pt]{elsarticle}

\usepackage{amssymb}
\usepackage{amsfonts}
\usepackage{amsmath}
\usepackage{theorem}
\usepackage{color}

\newtheorem{theorem}{Theorem}
{\theoremstyle{plain} \theorembodyfont{\rmfamily}
\newtheorem{definition}[theorem]{Definition}
\newtheorem{remark}[theorem]{Remark}
\newtheorem{example}[theorem]{Example}
}
{\theoremstyle{plain} \theorembodyfont{\rmfamily}

\newtheorem{corollary}[theorem]{Corollary}
\newtheorem{lemma}[theorem]{Lemma}
\newtheorem{proposition}[theorem]{Proposition}

}
{\theoremstyle{plain} \theorembodyfont{\rmfamily\bfseries}

}

\newenvironment{proof}[1][Proof]{\noindent\textbf{#1.} }{\ \rule{0.5em}{0.5em}}

\definecolor{pink}{rgb}{1,0.08,0.58}
\definecolor{orange}{rgb}{1,0.5,0}
\definecolor{purple}{rgb}{0.75002,0,1}
\definecolor{olive}{RGB}{85,107,47}
\definecolor{mygreen}{rgb}{0,0.6,0}

\begin{document}

\begin{frontmatter}

\title{\textbf{\textsf{Fuzzy ample spectrum contractions in (more general than) non-Archimedean fuzzy metric spaces}}}

\author[Antonio]{Antonio Francisco Rold\'{a}n L\'{o}pez de Hierro\corref{cor1}}   \ead{aroldan@ugr.es}
\author[Erdal,Erdal2]{Erdal Karap{\i}nar} \ead{erdalkarapinar@yahoo.com, karapinar@mail.cmuh.org.tw}
\author[Naseer]{Naseer Shahzad}   \ead{nshahzad@kau.edu.sa}

\cortext[cor1]{Corresponding author: Antonio Francisco Rold\'{a}n L\'{o}pez de Hierro, aroldan@ugr.es\\
\hspace*{5mm} 2010 Mathematics Subject Classification: 47H10, 47H09, 54H25, 46T99.}

\address[Antonio]{Department of Statistics and Operations Research, University of Granada, Granada, Spain.}
\address[Erdal]{Department of Medical Research, China Medical University, Taichung 40402, Taiwan.}
\address[Erdal2]{Department of Mathematics, \c{C}ankaya University, 06790, Etimesgut, Ankara, Turkey.}
\address[Naseer]{Department of Mathematics, Faculty of Science,  King Abdulaziz University, P.O.B. 80203, Jeddah 21589, Saudi Arabia.}

\date{October 22$^{th}$, 2020}

\begin{abstract}
Taking into account that Rold\'{a}n \emph{et al.}'s \emph{ample spectrum contractions} have managed to extend and unify more than ten distinct families of contractive mappings in the setting of metric spaces, in this manuscript we present a first study on how such concept can be implemented in the more general framework of fuzzy metric spaces in the sense of Kramosil and Mich\'{a}lek. We introduce two distinct approaches to the concept of \emph{fuzzy ample spectrum contractions} and we prove general results about existence and uniqueness of fixed points. The proposed notions enjoys the following advantages with respect to previous approaches: (1) they permit to develop fixed point theory in a very general fuzzy framework (for instance, the underlying fuzzy space is not necessarily complete); (2) the procedures that we employ are able to overcome the technical drawbacks raising in fuzzy metric spaces in the sense of Kramosil and Mich\'{a}lek (that do not appear in fuzzy metric spaces in the sense of George and Veeramani); (3) we introduce a novel property about sequences that are not Cauchy in a fuzzy space in order to consider a more general context than non-Archimedean fuzzy metric spaces; (4) the contractivity condition associated to \emph{fuzzy ample spectrum contractions} does not have to be expressed in separate variables; (5) such fuzzy contractions generalize some well-known families of fuzzy contractive operators such that the class of all Mihe\c{t}'s fuzzy $\psi$-contractions. As a consequence of these properties, this study gives a positive partial answer to a question posed by this author some years ago.
\end{abstract}


\end{frontmatter}

\section{\textbf{Introduction}}

\emph{Fixed point theory} has become, in recent years, into a very flourishing
branch of Nonlinear Analysis due, especially, to its ability to find solutions
of nonlinear equations. After Banach's pioneering definition, a multitude of
results have appeared that generalize his famous contractive mapping
principle. Some extensions generalize the contractive condition (see, e.g., \cite{AyKaRo}) and others
focus their efforts on the metric characteristics of the underlying space.

Fixed point theory has been successfully applied in many fields, especially in
Nonlinear Analysis, where it constitutes a major instrument in order to the existence and the uniqueness of solutions of several
classes of equations such as functional equations, matrix equations \cite{AsMo,Berzig2012}, integral equations \cite{AyNaSaYa,HaLoSa,SaHuShFaRa}, nonlinear systems \cite{AAKR}, etc.

In \cite{RoSh-ample1} Rold\'{a}n L\'{o}pez de Hierro and Shahzad introduced a
new family of contractions whose main characteristic was its capacity to
extend and unify several classes of well-known previous contractive mappings
in the setting of metric spaces: (1) Banach's contractions; (2) manageable
contractions \cite{DuKh}; (3) Farshid \emph{et al}'s $\mathcal{Z}%
$-contractions involving \emph{simulation functions} \cite{KhShRa,RoKaMa}; (4)
Geraghty's contractions \cite{Geraghty}; (5) Meir and Keeler's contractions
\cite{MeKe,Lim}; (6) $R$-contractions \cite{RoSh1}; (7) $(R,\mathcal{S}%
)$-contractions \cite{RoSh2}; (8) $(A,\mathcal{S})$-contractions
\cite{ShRoFa1}; (9) Samet et al's contractions \cite{SaVeVe}; (10) Shahzad
\emph{et al.}'s contractions \cite{ShKaRo}; (11) Wardowski's F-contractions
\cite{War}; etc. In a broad (or even philosophical) sense, the conditions that
define a Rold\'{a}n \emph{et al.}'s ample spectrum contraction should rather
be interpreted as the minimal set of properties that any contraction must
satisfy. Therefore, from our point of view, any other approach to contractive
mappings in the field of fixed point theory could take them into account.

One of the main settings in which more work has been done to extend contractive mappings in metric spaces is the context of fuzzy metric spaces (see, e.g., \cite{AlMi,GV,Grabiec,GrSa,Mi3,RoKaManro,FA-spaces}). These abstract spaces are able of providing an appropriate point of view for determining how similar or distinct two imprecise quantities are, which gives added value to the real analytical techniques that are capable of being extended to the fuzzy setting. In fact, the category of fuzzy metric spaces is so rich and broad that, as we shall see, the notion of fuzzy ample spectrum contraction does not directly come from its real corresponding counterpart. As a consequence of this great variety of examples, fuzzy metrics are increasingly appreciated in the scientific field in general due to the enormous popularity and importance of the results that are currently being obtained when working with fuzzy sets instead of classical sets, especially in Computation and Artificial Intelligence.

In this work we give a first introduction on how ample spectrum contractions
can be extended to fuzzy metric spaces. We will consider fuzzy metric spaces
in the sense of Kramosil and Mich\'{a}lek \cite{KrMi}, which are more general
than fuzzy metric spaces in the sense of George and Veeramani \cite{GV}.
Undoubtedly, the second kind of fuzzy metric spaces are easier to handle than
the first one because in the second case the metric only takes strictly
positive values. Kramosil and Mich\'{a}lek's fuzzy spaces are as general that
they include the case in which the distance between two points is infinite,
that is, its fuzzy metric is constantly zero (see \cite{MRR3}). This fact greatly complicates
the proofs that could be made of many fixed point theorems in the context of
George and Veeramani's fuzzy spaces (which has been a significant challenge).

In our study, we have tried to be as faithful as possible to the original
definition of ample spectrum contraction in metric spaces. This could mislead
the reader to the wrong idea that fuzzy ample spectrum contractions are
relatively similar to ample spectrum contractions in metric spaces. This is a
completely false statement: fuzzy ample spectrum contractions are so singular
that, for the moment, these authors have not been able to demonstrate that
each ample spectrum contraction in a metric space can be seen as a fuzzy ample
spectrum contractions in a fuzzy metric space (we have only proved it under
one additional condition). The difficulty appears when, in a fuzzy metric
space, we try to define a mapping that could extend the mapping $\varrho$ of
the ample spectrum contraction. In fact, as we commented above, the wide variety of
distinct classes of fuzzy metrics, and the way in which the fuzzy distance between
two points is represented by a distance distribution function, greatly complicates the
task of studying some possible relationship between ample spectrum contractions in
real metric spaces and in fuzzy metric spaces.

The presented fuzzy contractivity condition makes only use of two metric terms:
the fuzzy distance between two distinct points, $M(x,y,t)$, and the
fuzzy distance between their images, $M(Tx,Ty,t)$, under the self-mapping $T$.
Traditionally, these two terms have played separate roles in classical fuzzy and non-fuzzy
contractivity conditions: for instance, they use to appear in distinct sides of the inequality.
However, inspired by Farshid \emph{et al}'s $\mathcal{Z}$-contractions \cite{KhShRa,RoKaMa}
(based on \emph{simulation functions}), the contractivity condition associated to
\emph{fuzzy ample spectrum contractions} is not necessarily a constraint on separate variables.
As a consequence of this generalization, it is easy to check that the canonical examples
motivating the fixed point theory, that is, the Banach's contractive mappings, are particular cases of
\emph{fuzzy ample spectrum contractions}. Furthermore, we illustrate how Mihe\c{t} fuzzy
$\psi$-contractions \cite{Mi3} can also be seen as \emph{fuzzy ample spectrum contractions}, and why
this last class of fuzzy contractions are not the best choice to check that
Banach's contractions in metric spaces are fuzzy contractions.

In this paper we introduce two kind of fuzzy ample spectrum contractions. The
first one is more strict, but it allows to demonstrate some fixed point
results when the contractive condition does not satisfy additional
constraints. The second definition is more general but it forces us to handle
a concrete subfamily of contractivity conditions. We prove some distinct
results about existence and uniqueness of fixed points associated to each
class of fuzzy contractions. Our results generalize some previous fixed point
theorems introduced, on the one hand, by Mihe\c{t} and, on the other hand, by
Altun and Mihe\c{t}.

By the way, this study gives a positive partial answer to a question posed by Mihe\c{t} in \cite{Mi3} some years ago.
In such paper, the author wondered whether some fixed point theorems involving fuzzy $\psi$-contractions
could also hold in fuzzy metric spaces that did not necessarily satisfy the non-Archimedean property.
To face this challenge, we have introduced a novel assumption that we have called \emph{property} $\mathcal{NC}$ (\textquotedblleft non-Cauchy\textquotedblright).
Such condition establishes a very concrete behavior for asymptotically regular sequences that do not satisfy the Cauchy's condition:
in such cases, the fuzzy distance between two partial subsequences and their corresponding predecessors can be controlled
in terms of convergence. Obviously, our study is a proper generalization because each non-Archimedean fuzzy metric space satisfies such property
even if the associated triangular norm is only continuous at the boundary of the unit square (but not necessarily in the interior of the unit square).

\vspace*{3mm} 

\section{\textbf{Preliminaries}}

For an optimal understanding of this paper, we introduce here some basic
concepts and notations that could also be found in
\cite{G-book,RoSh-ample1,ShRoFa1}. Throughout this manuscript, let
$\mathbb{R}$ be the family of all real numbers, let $\mathbb{I}$ be the real
compact interval $\left[  0,1\right]  $, let $\mathbb{N}=\{1,2,3,\ldots\}$
denote the set of all positive integers and let $\mathbb{N}_{0}=\mathbb{N}%
\cup\{0\}$. Henceforth, $X$ will denote a non-empty set. A {\emph{binary
relation on }$X$} is a non-empty subset $\mathcal{S}$ of the Cartesian product
space $X\times X$. The notation $x\mathcal{S}y$ means that $(x,y)\in
\mathcal{S}$. We write $x\mathcal{S}^{\mathcal{\ast}}y$ when $x\mathcal{S}y$
and $x\neq y$. Hence $\mathcal{S}^{\mathcal{\ast}}$ is another binary relation
on $X$ (if it is non-empty). Two points $x$ and $y$ are $\mathcal{S}%
$\emph{-comparable} if $x\mathcal{S}y$ or $y\mathcal{S}x$. We say that
$\mathcal{S}$ is \emph{transitive} if we can deduce $x\mathcal{S}z$ from
$x\mathcal{S}y$ and $y\mathcal{S}z$. The \emph{trivial binary relation}
$\mathcal{S}_{X}$ on $X$ is defined by $x\mathcal{S}_{X}y$ for each $x,y\in
X$.

From now on, let $T:X\rightarrow X$ be a map from $X$ into itself. We say that
$T$ is {$\mathcal{S}$\emph{-nondecreasing}} if $Tx\mathcal{S}Ty$ for each
$x,y\in X$ such that $x\mathcal{S}y$. If a point $x\in X$ verifies $Tx=x $,
then $x$ is a \emph{fixed point of }$T$. We denote by $\operatorname*{Fix}(T)$
the set of all fixed points of $T$.

A sequence $\{x_{n}\}_{n\in\mathbb{N}_{0}}$ is called a {\emph{Picard sequence
of }$T$} \emph{based on }$x_{0}\in X$ if $x_{n+1}=Tx_{n}$ for all
$n\in\mathbb{N}_{0}$. Notice that, in such a case, $x_{n}=T^{n}x_{0}$ for each
$n\in\mathbb{N}_{0}$, where $\{T^{n}:X\rightarrow X\}_{n\in\mathbb{N}_{0}}$
are the \emph{iterates of }$T$ defined by $T^{0}=\,$identity, $T^{1}=T $ and
$T^{n+1}=T\circ T^{n}$ for all $n\geq2$.

Following \cite{RoShJNSA}, a sequence $\left\{  x_{n}\right\}  $ in $X$ is
\emph{infinite} if $x_{n}\neq x_{m}$ for all $n\neq m$, and $\left\{
x_{n}\right\}  $ is \emph{almost periodic} if there exist $n_{0}%
,N\in\mathbb{N}$ such that%
\[
x_{n_{0}+k+Np}=x_{n_{0}+k}\quad\text{for all }p\in\mathbb{N}\text{ and all
}k\in\left\{  0,1,2,\ldots,N-1\right\}  .
\]

\begin{proposition}
\label{K38 - 22 lem infinite or almost periodic}\textrm{(\cite{RoShJNSA},
Proposition 2.3)} Every Picard sequence is either infinite or almost periodic.
\end{proposition}

\subsection{\textbf{Ample spectrum contractions in metric spaces}}

We describe here the notion of \emph{ample spectrum contraction} in the
context of a metric space. Such concept involves a key kind of sequences of
real numbers that must be highlighted. Let $(X,d)$ be a metric space, let
$\mathcal{S}$ be a binary relation on $X$, let $T:X\rightarrow X$ be a
self-mapping and let $\varrho:A\times A\rightarrow\mathbb{R}$ be a function
where $A\subseteq\mathbb{R}$ is a non-empty subset of real numbers.

\begin{definition}
\label{K38 - 02 def TS sequence in MS}\textrm{(\cite{RoSh-ample1}, Definition
3)} Let $\{a_{n}\}$ and $\{b_{n}\}$ be two sequences of real numbers. We say
that $\{\left(  a_{n},b_{n}\right)  \}$ is a {$(T,\mathcal{S}^{\ast}%
)$\emph{-sequence}} if there exist two sequences $\{x_{n}\},\{y_{n}\}\subseteq
X$ such that%
\[
x_{n}\mathcal{S}^{\ast}y_{n},\quad Tx_{n}\mathcal{S}^{\ast}Ty_{n},\quad
a_{n}=d(Tx_{n},Ty_{n})\quad\text{and}\quad b_{n}=d(x_{n},y_{n})\quad\text{for
all }n\in\mathbb{N}_{0}.
\]

\end{definition}

The previous class of sequences plays a crucial role in the third condition of
the following definition.

\begin{definition}
\label{K38 - 01 def ample spectrum contraction MS}\textrm{(\cite{RoSh-ample1},
Definition 4)} We will say that $T:X\rightarrow X$ is an \emph{ample spectrum
contraction} w.r.t. $\mathcal{S}$ and $\varrho$ if the following four
conditions are fulfilled.

\begin{description}
\item[$\left(  \mathcal{B}_{1}\right)  $] $A$ is nonempty and $\left\{
\,d\left(  x,y\right)  \in\left[  0,\infty\right)  :x,y\in X,~x\mathcal{S}%
^{\ast}y\,\right\}  \subseteq A$.

\item[$\left(  \mathcal{B}_{2}\right)  $] If $\{x_{n}\}\subseteq X$ is a
Picard $\mathcal{S}$-nondecreasing sequence of $T$ such that
\[
x_{n}\neq x_{n+1}\quad\text{and}\quad\varrho\left(  d\left(  x_{n+1}%
,x_{n+2}\right)  ,d\left(  x_{n},x_{n+1}\right)  \right)  \geq0\quad\text{for
all }n\in\mathbb{N}_{0},
\]
then $\{d\left(  x_{n},x_{n+1}\right)  \}\rightarrow0$.

\item[$\left(  \mathcal{B}_{3}\right)  $] If $\{\left(  a_{n},b_{n}\right)
\}\subseteq A\times A$ is a $(T,\mathcal{S}^{\ast})$-sequence such that
$\{a_{n}\}$ and $\{b_{n}\}$ converge to the same limit $L\geq0$ and verifying
that $L<a_{n}$ and $\varrho(a_{n},b_{n})\geq0$ for all $n\in\mathbb{N}_{0}$,
then $L=0$.

\item[$\left(  \mathcal{B}_{4}\right)  $] {$\varrho\left(
d(Tx,Ty),d(x,y)\right)  \geq0\quad$for all $x,y\in X\quad$such that
$x\mathcal{S}^{\ast}y$ and $Tx\mathcal{S}^{\ast}Ty$.}
\end{description}
\end{definition}

In some cases, we will also consider the following two auxiliary properties.

\begin{description}
\item[$\left(  \mathcal{B}_{2}^{\prime}\right)  $] If $x_{1},x_{2}\in X$ are
two points such that%
\[
T^{n}x_{1}\mathcal{S}^{\ast}T^{n}x_{2}\quad\text{and}\quad{\varrho({d}\left(
T^{n+1}x_{1},T^{n+1}x_{2}\right)  ,{d}\left(  T^{n}x_{1},T^{n}x_{2}\right)
)\geq0}\quad\text{for all }n\in\mathbb{N}_{0},
\]
then $\{{d}\left(  T^{n}x_{1},T^{n}x_{2}\right)  \}\rightarrow0$.

\item[$\left(  \mathcal{B}_{5}\right)  $] If $\{\left(  a_{n},b_{n}\right)
\}$ is a $(T,\mathcal{S}^{\ast})$-sequence such that $\{b_{n}\}\rightarrow0$
and $\varrho(a_{n},b_{n})\geq0$ for all $n\in\mathbb{N}_{0}$, then
$\{a_{n}\}\rightarrow0$.
\end{description}

Rold\'{a}n L\'{o}pez de Hierro and Shahzad demonstrated that, under very weak
conditions, these contractions have a fixed point and, if we assume other
constraints, then such fixed point is unique (see \cite{RoSh-ample1}). After
that, the same authors and Karap\i nar were able to extend their study to
Branciari distance spaces.

\subsection{\textbf{Triangular norms}}

A \emph{triangular norm} \cite{ScSk}\ (for short, a \emph{t-norm}) is a
function $\ast:\mathbb{I}\times\mathbb{I}\rightarrow\mathbb{I}$ satisfying the
following properties: associativity, commutativity, non-decreasing on each
argument, has $1$ as unity (that is, $t\ast1=t$ for all $t\in\mathbb{I}$). It
is usual that authors consider continuous t-norms on their studies. A t-norm
$\ast$ is \emph{positive} if $t\ast s>0$ for all $t,s\in\left(  0,1\right]  $.
Given two t-norms $\ast$ and $\ast^{\prime}$, we will write $\ast\leq
\ast^{\prime}$ when $t\ast s\leq t\ast^{\prime}s$ for all $t,s\in\mathbb{I}$.
Examples of t-norms are the following ones:%
\[%
\begin{tabular}
[c]{ll}%
$\text{Product }\ast_{P}:$ & $t\ast_{P}s=t\,s$\\
$\text{\L ukasiewicz }\ast_{L}:$ & $t\ast_{L}s=\max\{0,t+s-1\}$\\
$\text{Minimum }\ast_{m}:$ & $t\ast_{m}s=\min\{t,s\}$\\
$\text{Drastic }\ast_{D}:$ & $t\ast_{D}s=\left\{
\begin{tabular}
[c]{ll}%
$0,$ & if $t<1$ and $s<1,$\\
$\min\{t,s\},$ & if $t=1$ or $s=1.$%
\end{tabular}
\right.  $%
\end{tabular}
\]
If $\ast$ is an arbitrary t-norm, then $\ast_{D}\leq\ast\leq\ast_{m}$, that
is, the drastic t-norm is the absolute minimum and the minimum t-norm is the
absolute maximum among the family of all t-norms (see \cite{t-norm}).

\begin{definition}
We will say that a t-norm $\ast$ is \emph{continuous at the }$1$%
\emph{-boundary} if it is continuous at each point of the type $\left(
1,s\right)  $ where $s\in\mathbb{I}$ (that is, if $\{t_{n}\}\rightarrow1$ and
$\{s_{n}\}\rightarrow s$, then $\{t_{n}\ast s_{n}\}\rightarrow1\ast s=s$).
\end{definition}

Obviously, each continuous t-norm is continuous at the $1$-boundary.

\begin{proposition}
\label{K38 - 18 propo cancellation}Let $\{a_{n}\},\{b_{n}\},\{c_{n}%
\},\{d_{n}\},\{e_{n}\}\subseteq\mathbb{I}$ be five sequences and let
$L\in\mathbb{I}$ be a number such that $\{a_{n}\}\rightarrow L$,
$\{b_{n}\}\rightarrow1$, $\{d_{n}\}\rightarrow1$ and $\{e_{n}\}\rightarrow L$.
Suppose that $\ast$ is a continuous at the $1$-boundary $t$-norm and that
\[
a_{n}\geq b_{n}\ast c_{n}\ast d_{n}\geq e_{n}\quad\text{for all }%
n\in\mathbb{N}.
\]
Then $\{c_{n}\}$ converges to $L$.
\end{proposition}

\begin{proof}
As $\{c_{n}\}\subseteq\mathbb{I}=\left[  0,1\right]  $ is bounded, then it has
a convergent partial subsequence. Let $\{c_{\sigma(n)}\}$ be an arbitrary
convergent partial subsequence of $\{c_{n}\}$ and let $L^{\prime}%
=\lim_{n\rightarrow\infty}c_{\sigma(n)}$. Since $\ast$ is continuous at the
point $\left(  1,L^{\prime}\right)  $, $\{b_{\sigma(n)}\}\rightarrow1$ and
$\{c_{\sigma(n)}\}\rightarrow L^{\prime}$, then $\{b_{\sigma(n)}\ast
c_{\sigma(n)}\}\rightarrow1\ast L^{\prime}=L^{\prime}$, and as $\{d_{\sigma
(n)}\}\rightarrow1$, then $\{b_{\sigma(n)}\ast c_{\sigma(n)}\ast d_{\sigma
(n)}\}\rightarrow L^{\prime}\ast1=L^{\prime}$. Furthermore, taking into
account that $a_{n}\geq b_{n}\ast c_{n}\ast d_{n}\geq e_{n}$ for all
$n\in\mathbb{N}$, we deduce that%
\[
L=\lim_{n\rightarrow\infty}a_{\sigma(n)}\geq\lim_{n\rightarrow\infty}\left(
b_{\sigma(n)}\ast c_{\sigma(n)}\ast d_{\sigma(n)}\right)  \geq\lim
_{n\rightarrow\infty}e_{\sigma(n)}=L.
\]
Hence $L^{\prime}=\lim_{n\rightarrow\infty}\left(  b_{\sigma(n)}\ast
c_{\sigma(n)}\ast d_{\sigma(n)}\right)  =L$. This proves that any convergent
partial subsequence of $\{c_{n}\}$ converges to $L$. Next we consider the
limit inferior and the limit superior of $\{c_{n}\}$. The previous argument
shows that%
\[
L=\liminf_{n\rightarrow\infty}c_{n}\leq\limsup_{n\rightarrow\infty}c_{n}=L.
\]
As the limit inferior and the limit superior of $\{c_{n}\}$ are equal to $L$,
the sequence $\{c_{n}\}$ is convergent and its limit is $L$.
\end{proof}

\begin{example}
The cancelation property showed in Proposition
\ref{K38 - 18 propo cancellation} is not satisfied by all t-norms. For
instance, let $\ast_{D}$ be the drastic t-norm (which is not continuous at any
point $(1,s)$ or $\left(  s,1\right)  $ of the $1$-boundary when $s>0$). Let
$L=0$ and let $\{a_{n}\},\{b_{n}\},\{c_{n}\},\{d_{n}\},\{e_{n}\}\subseteq
\mathbb{I}$ be the sequences on $\mathbb{I}$ given, for all $n\in\mathbb{N}$,
by:%
\[
a_{n}=e_{n}=0,\quad b_{n}=d_{n}=1-\frac{1}{n},\quad c_{n}=\frac{1}{2}.
\]
Then $b_{n}\ast_{D}c_{n}=b_{n}\ast_{D}c_{n}\ast_{D}d_{n}=a_{n}=e_{n}=L=0$ for
all $n\in\mathbb{N}$, $\{b_{n}\}\rightarrow1$ and $\{d_{n}\}\rightarrow1$.
However, $\{c_{n}\}$ does not converge to $L=0$. In fact, it has not any partial
subsequence converging to $L=0$.
\end{example}

\subsection{\textbf{Fuzzy metric spaces}}

In this subsection we introduce two distinct notions of \emph{fuzzy metric
space} that represent natural extensions of the concept of \emph{metric space}
to a setting in which some uncertainty or imprecision can be considered when
determining the distance between two points.

\begin{definition}
\label{definition KM-space}\textrm{(cf. Kramosil and Mich\'{a}lek
\cite{KrMi})} A triplet $(X,M,\ast$) is called a \emph{fuzzy metric space in
the sense of Kramosil and Mich\'{a}lek} (briefly, a \emph{KM-FMS}) if $X$ is
an arbitrary non-empty set, $\ast$ is a $t$-norm and $M:X\times X\times\left[
0,\infty\right)  \rightarrow\mathbb{I} $ is a fuzzy set satisfying the
following conditions, for each $x,y,z\in X$, and $t,s\geq0$:

\begin{description}
\item[(KM-1)] $M(x,y,0)=0$;

\item[(KM-2)] $M(x,y,t)=1$ for all $t>0$ if, and only if, $x=y$;

\item[(KM-3)] $M(x,y,t)=M(y,x,t)$;

\item[(KM-4)] $M(x,z,t+s)\geq M(x,y,t)\ast M(y,z,s)$;

\item[(KM-5)] $M(x,y,\cdot):\left[  0,\infty\right)  \rightarrow\left[
0,1\right]  $ is left-continuous.
\end{description}
\end{definition}

The value $M(x,y,t)$ can be interpreted of as the degree of nearness between
$x$ and $y$ compared to $t$. On their original definition, Kramosil and
Mich\'{a}lek did not assume the continuity of the t-norm $\ast$. However, in
later studies in KM-FMS, it is very usual to suppose that $\ast$ is continuous
(see, for instance, \cite{Mi3}).

The following one is the canonical way in which a metric space can be seen as
a KM-FMS.

\begin{example}
\label{K38 - 24 ex canonical metric FMS}Each metric space $\left(  X,d\right)
$ can be seen as a KM-FMS $(X,M^{d},\ast)$, where $\ast$ is any t-norm, by
defining $M:X\times X\times\left[  0,\infty\right)  \rightarrow\mathbb{I}$ as:%
\[
M^{d}\left(  x,y,t\right)  =\left\{
\begin{tabular}
[c]{ll}%
$0,$ & if $t=0,$\\
$\dfrac{t}{t+d\left(  x,y\right)  },$ & if $t>0.$%
\end{tabular}
\right.
\]
Notice that $0<M^{d}\left(  x,y,t\right)  <1$ for all $t>0$ and all $x,y\in X$
such that $x\neq y$. Furthermore, $\lim_{t\rightarrow\infty}M^{d}\left(
x,y,t\right)  =1$ for all $x,y\in X$. More properties of these spaces are
given in Proposition \ref{K38 - 44 propo Md non-Archimedean}.
\end{example}

\begin{remark}
\label{K38 - 32 rem Km to GV}Definition \ref{definition KM-space} is as
general that such class of fuzzy spaces can verify that%
\begin{equation}
M(x,y,t)=0\text{\quad for all }t>0\text{\quad when }x\neq y.
\label{K38 - 26 prop bad condition}%
\end{equation}
In the context of \emph{extended metric spaces} (whose metrics take values in
the extended interval $\left[  0,\infty\right]  $, including $\infty$; see
\cite{MRR3,FA-spaces}), this property can be interpreted by saying that the
\textquotedblleft distance\textquotedblright\ between the points $x$ and $y$
is infinite. Property (\ref{K38 - 26 prop bad condition}) is usually
unsuitable in the setting of fixed point theory because it often spoils the
arguments given in the proofs. This drawback is sometimes overcame by assuming
that the initial condition is a point $x_{0}\in X$ such that $M(x_{0}%
,Tx_{0},t)>0$ for all $t>0$ because the contractivity condition helps to prove
that all the points of the Picard sequence based on $x_{0}$ satisfy the same
condition. In order to hereby such characteristic, we would need to assume
additional constraints on the auxiliary functions we will employ to introduce
the announced fuzzy ample spectrum contractions. Such constraints would not
appear if we would have decided to restrict out study to fuzzy metric spaces
in the sense of George and Veeramani, which were introduced in order to
consider a Hausdorff topology on the corresponding fuzzy spaces and to prove a
version of the Baire theorem.
\end{remark}

\begin{definition}
\textrm{(cf. George and Veeramani \cite{GV})} A triplet $(X,M,\ast$) is called
a \emph{fuzzy metric space in the sense of George and Veeramani} (briefly,
a\emph{\ GV-FMS}) if $X$ is an arbitrary non-empty set, $\ast$ is a continuous
$t$-norm and $M:X\times X\times\left(  0,\infty\right)  \rightarrow\mathbb{I}$
is a fuzzy set satisfying the following conditions, for each $x,y,z\in X$, and
$t,s>0$:

\begin{description}
\item[(GV-1)] $M(x,y,t)>0$;

\item[(GV-2)] $M(x,y,t)=1\text{ for all }t>0$ if, and only if, $x=y$;

\item[(GV-3)] $M(x,y,t)=M(y,x,t)$;

\item[(GV-4)] $M(x,z,t+s)\geq M(x,y,t)\ast M(y,z,s)$;

\item[(GV-5)] $M(x,y,\cdot):\left(  0,\infty\right)  \rightarrow\left[
0,1\right]  $ is a continuous function.
\end{description}
\end{definition}

\begin{lemma}
Each GV-FMS is, in fact, a KM-FMS by extending $M$ to $t=0$ as $M(x,y,0)=0$
for all $x,y\in X$.
\end{lemma}

\begin{lemma}
\textrm{(cf. Grabiec \cite{Grabiec})} If $\left(  X,M,\ast\right)  $ is a
KM-FMS (respectively, a GV-FMS) and $x,y\in X$, then each function
$M(x,y,\cdot)$ is nondecreasing on $\left[  0,\infty\right)  $ ((respectively,
on $\left(  0,\infty\right)  $).
\end{lemma}

\begin{proposition}
\label{K40 21 propo either infinite or almost-constant fuzzy}\textrm{(cf.
\cite{RoKaFu}, Proposition 2)} Let $\left\{  x_{n}\right\}  $ be a Picard
sequence in a KM-FMS $(X,M,\ast)$ such that $\{M(x_{n},x_{n+1},t)\}\rightarrow
1$ for all $t>0$. If there are $n_{0},m_{0}\in\mathbb{N}$ such that
$n_{0}<m_{0}$ and $x_{n_{0}}=x_{m_{0}}$, then there is $\ell_{0}\in\mathbb{N}$
and $z\in X$ such that $x_{n}=z$ for all $n\geq\ell_{0}$ (that is, $\left\{
x_{n}\right\}  $ is constant from a term onwards). In such a case, $z$ is a
fixed point of the self-mapping for which $\left\{  x_{n}\right\}  $ is a
Picard sequence.
\end{proposition}

\begin{proof}
Let $T:X\rightarrow X$ be a mapping for which $\left\{  x_{n}\right\}  $ is a
Picard sequence, that is, $x_{n+1}=Tx_{n}$ for all $n\in\mathbb{N}$. The set
\[
\Omega=\left\{  \,k\in\mathbb{N}:\exists~n_{0}\in\mathbb{N}\text{ such that
}x_{n_{0}}=x_{n_{0}+k}\,\right\}
\]
is non-empty because $m_{0}-n_{0}\in\Omega$, so it has a minimum $k_{0}%
=\min\Omega$. Then $k_{0}\geq1$ and there is $n_{0}\in\mathbb{N}$ such that
$x_{n_{0}}=u_{n_{0}+k_{0}}$. As $\left\{  x_{n}\right\}  $ is not infinite,
then it must be almost periodic by Proposition
\ref{K38 - 22 lem infinite or almost periodic}. In fact, it is easy to check,
by induction on $p$, that:%
\begin{equation}
x_{n_{0}+r+pk_{0}}=x_{n_{0}+r}\quad\text{for all }p\in\mathbb{N}\text{ and all
}r\in\left\{  0,1,2,\ldots,k_{0}-1\right\}  . \label{K40 20 prop}%
\end{equation}

If $k_{0}=1$, then $x_{n_{0}}=x_{n_{0}+1}$. Similarly $x_{n_{0}+2}%
=Tx_{n_{0}+1}=Tx_{n_{0}}=x_{n_{0}+1}=x_{n_{0}}$. By induction, $x_{n_{0}%
+r}=x_{n_{0}}$ for all $r\geq0$, which is precisely the conclusion. Next we
are going to prove that the case $k_{0}\geq2$ leads to a contradiction.

Assume that $k_{0}\geq2$. Then each two terms in the set $\{\,x_{n_{0}%
},x_{n_{0}+1},x_{n_{0}+2},\ldots,x_{n_{0}+k_{0}-1}\,\}$ are distinct, that is,
$x_{n_{0}+i}\neq x_{n_{0}+j}$ for all $0\leq i<j\leq k_{0}-1$ (on the contrary
case, $k_{0}$ is not the minimum of $\Omega$). Since $x_{n_{0}+i}\neq
x_{n_{0}+i+1}$ for all $i\in\left\{  0,1,2,\ldots,k_{0}-1\right\}  $, then
there is $s_{i}>0$ such that $M(x_{n_{0}+i},x_{n_{0}+i+1},s_{i})<1$ for all
$i\in\left\{  0,1,2,\ldots,k_{0}-1\right\}  $. Since each $M(x_{n_{0}%
+i},x_{n_{0}+i+1},\cdot)$ is a non-decreasing function, if $t_{0}%
=\min(\{\,s_{i}:0\leq i\leq k_{0}-1\,\})>0$, then%
\[
M(x_{n_{0}+i},x_{n_{0}+i+1},t_{0})<1\quad\text{for all }i\in\left\{
0,1,2,\ldots,k_{0}-1\right\}  .
\]
Let define%
\[
\delta_{0}=\max\left(  \left\{  \,M(x_{n_{0}+i},x_{n_{0}+i+1},t_{0}):0\leq
i\leq k_{0}-1\,\right\}  \right)  \in\left[  0,1\right)  .
\]
Then $\delta_{0}<1$. Since $\lim_{n\rightarrow\infty}M(x_{n},x_{n+1},t_{0})=1
$, there is $r_{0}\in\mathbb{N}$ such that $r_{0}\geq n_{0}$ and $M(x_{r_{0}%
},x_{r_{0}+1},t_{0})>\delta_{0}$. Let $i_{0}\in\left\{  0,1,2,\ldots
,k_{0}-1\right\}  $ be the unique integer number such that the non-negative
integer numbers $r_{0}-n_{0}$ and $i_{0}$ are congruent modulo $k_{0}$, that
is, $i_{0}$ is the rest of the integer division of $r_{0}-n_{0}$ over $k_{0}$.
Hence there is a unique integer $p\geq0$ such that $\left(  r_{0}%
-n_{0}\right)  -i_{0}=pk_{0}$. Since $r_{0}=n_{0}+i_{0}+pk_{0}$, property
(\ref{K40 20 prop}) guarantees that%
\[
x_{r_{0}}=x_{n_{0}+i_{0}+pk_{0}}=x_{n_{0}+i_{0}},
\]
where $n_{0}+i_{0}\in\left\{  n_{0},n_{0}+1,n_{0}+2,\ldots,n_{0}%
+k_{0}-1\right\}  $. As a consequence:%
\[
\delta_{0}=\max\left(  \left\{  \,M(x_{n_{0}+i},x_{n_{0}+i+1},t_{0}):0\leq
i\leq k_{0}-1\,\right\}  \right)  \geq M(x_{n_{0}+i_{0}},x_{n_{0}+i_{0}%
+1},t_{0})=M(x_{r_{0}},x_{r_{0}+1},t_{0})>\delta_{0},
\]
which is a contradiction.
\end{proof}

\subsection{\textbf{Non-Arquimedean fuzzy metric spaces}}

A KM-FMS $(X,M,\ast)$ is said to be \emph{non-Archimedean} \cite{Istr}\ if%
\begin{equation}
M(x,z,t)\geq M(x,y,t)\ast M(y,z,t)\quad\text{for all }x,y,z\in X\text{ and all
}t>0. \label{K38 - 27 prop non-Arquimedean}%
\end{equation}
This property is equivalent to:%
\[
M(x,z,\max\{t,s\})\geq M(x,y,t)\ast M(y,z,s)\quad\text{for all }x,y,z\in X\text{
and all }t,s>0.
\]
Notice that this property depends on both the fuzzy metric $M$ and the t-norm
$\ast$. Furthermore, the non-Archimedean property
(\ref{K38 - 27 prop non-Arquimedean}) implies the triangle inequality (KM-4),
so each non-Arquimedean fuzzy metric space is a KM-FMS.

\begin{proposition}
\label{K38 - 44 propo Md non-Archimedean}Given a metric space $\left(
X,d\right)  $, let $(X,M^{d})$ be the canonical way to see $\left(
X,d\right)  $ as a KM-FMS (recall Example
\ref{K38 - 24 ex canonical metric FMS}). Then the following properties are fulfilled.

\begin{enumerate}
\item \label{K38 - 44 propo Md non-Archimedean, item 1}$(X,M^{d})$ is a KM-FMS
(and also a GV-FMS) under any t-norm $\ast$.

\item \label{K38 - 44 propo Md non-Archimedean, item 2}If $\ast$ is a t-norm
such that $\ast\leq\ast_{P}$, then $(X,M,\ast)$ is a non-Archimedean KM-FMS.

\item \label{K38 - 44 propo Md non-Archimedean, item 3}The metric space
$(X,d)$ satisfies $d\left(  x,z\right)  \leq\max\{d\left(  x,y\right)
,d(y,z)\}$ for all $x,y,z\in X$ if, and only if, $(X,M^{d},\ast_{m})$ is a
non-Archimedean KM-FMS.
\end{enumerate}
\end{proposition}

Notice that the \emph{discrete metric} on $X$, defined by $d\left(
x,y\right)  =0$ if $x=y$ and $d\left(  x,y\right)  =1$ if $x\neq y$, is an
example of metric satisfying the property involved in item
\ref{K38 - 44 propo Md non-Archimedean, item 3} of Proposition
\ref{K38 - 44 propo Md non-Archimedean}.

\begin{example}
\textrm{(Altun and Mihe\c{t} \cite{AlMi}, Example 1.3)} Let $\left(
X,d\right)  $ be a metric space and let $\vartheta$ be a nondecreasing and
continuous function from $\left(  0,\infty\right)  $ into $\left(  0,1\right)
$ such that $\lim_{t\rightarrow\infty}\vartheta(t)=1$. Let $\ast$ be a t-norm
such that $\ast\leq\ast_{P}$. For each $x,y\in X$ and all $t\in\left(
0,\infty\right)  $ , define%
\[
M(x,y,t)=\left[  \vartheta(t)\right]  ^{d\left(  x,y\right)  }.
\]
Then $\left(  X,M,\ast\right)  $ is a non-Archimedean KM-FMS.
\end{example}

\begin{example}
A KM-FMS is called \emph{stationary} when each function $t\mapsto M(x,y,t)$
does not depend on $t$, that is, it only depends in the points $x$ and $y$.
For instance, if $X=\left(  0,\infty\right)  $ and $M$ is defined by%
\[
M(x,y,t)=\frac{\min\{x,y\}}{\max\{x,y\}}\quad\text{for all }x,y\in\left(
0,\infty\right)  \text{ and all }t>0,
\]
then $(X,M,\ast)$ is a stationary KM-FMS. In fact, it is non-Archimedean.
\end{example}

\section{\textbf{Fuzzy spaces}}

There are many properties which are defined for mathematical objects in a
fuzzy metric space that, in fact, only depend on the fuzzy metric $M$, but not
on the t-norm $\ast$. For instance, the notions of convergency and Cauchy
sequence. Therefore, it is worth-noting to introduce such notions when $M$ is
an arbitrary function.

\begin{definition}
A \emph{fuzzy space} is a pair $(X,M)$ where $X$ is a non-empty set and $M$ is
a fuzzy set on $X\times X\times\left[  0,\infty\right)  $, that is, a mapping
$M:X\times X\times\left[  0,\infty\right)  \rightarrow\mathbb{I}$ (notice that
no additional conditions are assumed for $M$). In a fuzzy space $(X,M)$, we
say that a sequence $\{x_{n}\}\subseteq X$ is:

\begin{itemize}
\item $M$\emph{-Cauchy} if for all $\varepsilon\in\left(  0,1\right)  $ and
all $t>0$ there is $n_{0}\in\mathbb{N}$ such that $M\left(  x_{n}%
,x_{m},t\right)  >1-\varepsilon$ for all $n,m\geq n_{0}$;

\item $M$\emph{-convergent to }$x\in X$ if for all $\varepsilon\in\left(
0,1\right)  $ and all $t>0$ there is $n_{0}\in\mathbb{N}$ such that $M\left(
x_{n},x,t\right)  >1-\varepsilon$ for all $n\geq n_{0}$ (in such a case, we
write $\{x_{n}\}\rightarrow x$).
\end{itemize}

We say that the fuzzy space $(X,M)$ is \emph{complete} (or $X$ is
$M$\emph{-complete}) if each $M$-Cauchy sequence in $X$ is $M$-convergent to a
point of $X$.
\end{definition}

\begin{proposition}
\label{K38 - 49 propo uniqueness of limit}The limit of an $M$-convergent
sequence in a KM-FMS whose t-norm is continuous at the $1$-boundary is unique.
\end{proposition}

\begin{remark}
\begin{enumerate}
\item Using the previous definitions, it is possible to prove that a sequence
$\{x_{n}\}$ in a metric space $(X,d)$ is Cauchy (respectively, convergent to
$x\in X$) if, and only if, $\{x_{n}\}$ is $M^{d}$-Cauchy (respectively,
$M^{d}$-convergent to $x\in X$) in $(X,M^{d}$).

\item Notice that a sequence $\{x_{n}\}$ is $M$-convergent to $x\in X$ if, and
only if, $\lim_{n\rightarrow\infty}M(x_{n},x,t)=1$ for all $t>0$.
\end{enumerate}
\end{remark}

It is clear that one of the possible drawbacks of the previous definition is
the fact that, if $M$ does not satisfy additional properties such as
(KM-1)-(KM-5), then the limit of an $M$-convergent sequence has not to be
unique, or an $M$-convergent sequence needs not be $M$-Cauchy. However, it is
a good idea to consider general fuzzy spaces because such spaces allows us to
reflex about what properties depend on $M$ and $\ast$, and what other
properties only depend on $M$. In the second case, the next conditions will be
of importance in what follows. Notice that the following notions have sense
even if the limits of $M$-convergent sequences are not unique.

\begin{definition}
Let $\mathcal{S}$ be a binary relation on a fuzzy space $\left(  X,M\right)
$, let $Y\subseteq X$ be a nonempty subset, let $\{x_{n}\}$ be a sequence in
$X$ and let $T:X\rightarrow X$ be a self-mapping. We say that:

\begin{itemize}
\item $T$ is $\mathcal{S}$\emph{-nondecreasing-continuous} if $\{Tx_{n}%
\}\rightarrow Tz$ for all $\mathcal{S}$-nondecreasing sequence $\{x_{n}%
\}\subseteq X$ such that $\{x_{n}\}\rightarrow z\in X$;

\item $T$ is $\mathcal{S}$\emph{-strictly-increasing-continuous} if
$\{Tx_{n}\}\rightarrow Tz$ for all $\mathcal{S}$-strictly-increasing sequence
$\{x_{n}\}\subseteq X$ such that $\{x_{n}\}\rightarrow z\in X$;

\item $Y$ is $(\mathcal{S},M)$\emph{-strictly-increasing-complete} if every
$\mathcal{S}$-strictly-increasing and $M$-Cauchy sequence $\{y_{n}\}\subseteq
Y$ is $M$-convergent to a point of $Y$;

\item $Y$ is $(\mathcal{S},M)$\emph{-strictly-increasing-precomplete} if there
exists a set $Z$ such that $Y\subseteq Z\subseteq X$ and $Z$ is $(\mathcal{S}%
,M)$-strictly-increasing-complete;

\item $\left(  X,M\right)  $ is $\mathcal{S}$%
\emph{-strictly-increasing-regular} if, for all $\mathcal{S}$%
-strictly-increasing sequence $\{x_{n}\}\subseteq X$ such that $\{x_{n}%
\}\rightarrow z\in X$, it follows that $x_{n}\mathcal{S}z$ for all
$n\in\mathbb{N}$;

\item $\left(  X,M\right)  $ is \emph{metrically-}$\mathcal{S}$%
\emph{-strictly-increasing-regular} if, for all $\mathcal{S}$%
-strictly-increasing sequence $\{x_{n}\}\subseteq X$ such that $\{x_{n}%
\}\rightarrow z\in X$ and%
\[
M\left(  x_{n},x_{n+1},t\right)  >0\quad\text{for all }n\in\mathbb{N}\text{
and all }t>0,
\]
it follows that $x_{n}\mathcal{S}z$ and $M(x_{n},z,t)>0$ for all
$n\in\mathbb{N}$ and all $t>0$.
\end{itemize}
\end{definition}

The reader can notice that we will only use the previous notions when we are
working with infinite Picard sequences, so they could be refined in prospect work.

\section{\textbf{The property $\mathcal{NC}$}}

In a fuzzy metric space $(X,M,\ast)$, if $\{x_{n}\}$ is not an $M$-Cauchy
sequence, then there are $\varepsilon_{0}\in\left(  0,1\right)  $ and
$t_{0}>0$ such that, for all $k\in\mathbb{N}$, there are natural numbers
$m\left(  k\right)  ,n\left(  k\right)  \geq k$ such that $M\left(
x_{n(k)},x_{m(k)},t_{0}\right)  \leq1-\varepsilon_{0}$. Equivalently, there
are two partial subsequences $\{x_{n(k)}\}_{k\in\mathbb{N}}$ and
$\{x_{m(k)}\}_{k\in\mathbb{N}}$ of $\{x_{n}\}$ such that $k<n\left(  k\right)
<m\left(  k\right)  <n\left(  k+1\right)  $ and%
\[
1-\varepsilon_{0}\geq M\left(  x_{n(k)},x_{m(k)},t_{0}\right)  \quad\text{for
all }k\in\mathbb{N}.
\]
Associated to each number $n\left(  k\right)  $, if we take $m\left(
k\right)  $ is the least integer satisfying the previous property, then we can
also suppose that%
\[
M\left(  x_{n\left(  k\right)  },x_{m\left(  k\right)  -1},t_{0}\right)
>1-\varepsilon_{0}\geq M\left(  x_{n\left(  k\right)  },x_{m\left(  k\right)
},t_{0}\right)  \quad\text{for all }k\in\mathbb{N}.
\]
From these inequalities, in a general fuzzy metric space, it is difficult to
go further, even if we try to use the properties of the t-norm $\ast$. Next,
we are going to introduce a new condition on the fuzzy space in order to
guarantee that the sequences $\{M\left(  x_{n\left(  k\right)  },x_{m\left(
k\right)  },t_{0}\right)  \}_{k\in\mathbb{N}}$ and $\{M\left(  x_{n\left(
k\right)  -1},x_{m\left(  k\right)  -1},t_{0}\right)  \}_{k\in\mathbb{N}}$
satisfy additional properties. This will be of great help when we handle such
sequences in the proofs of fixed point theorems. Immediately after, we show
that non-Archimedean fuzzy KM-metric spaces satisfy this new condition (which
is not trivial at all).

\begin{definition}
We will say that a fuzzy space $(X,M)$ satisfies the \emph{property
}$\mathcal{NC}$ (\textquotedblleft\emph{not Cauchy}\textquotedblright) if for
each sequence $\{x_{n}\}\subseteq X$ which is not $M$-Cauchy and verifies
$\lim_{n\rightarrow\infty}M\left(  x_{n},x_{n+1},t\right)  =1$ for all $t>0$,
there are $\varepsilon_{0}\in\left(  0,1\right)  $, $t_{0}>0$ and two partial
subsequences $\{x_{n(k)}\}_{k\in\mathbb{N}}$ and $\{x_{m(k)}\}_{k\in
\mathbb{N}}$ of $\{x_{n}\}$ such that, for all $k\in\mathbb{N}$,%
\begin{align*}
&  k<n\left(  k\right)  <m\left(  k\right)  <n\left(  k+1\right)
\qquad\text{and}\\
&  M\left(  x_{n\left(  k\right)  },x_{m\left(  k\right)  -1},t_{0}\right)
>1-\varepsilon_{0}\geq M\left(  x_{n\left(  k\right)  },x_{m\left(  k\right)
},t_{0}\right)  ,
\end{align*}
and also%
\[
\lim_{k\rightarrow\infty}M\left(  x_{n\left(  k\right)  },x_{m\left(
k\right)  },t_{0}\right)  =\lim_{k\rightarrow\infty}M\left(  x_{n\left(
k\right)  -1},x_{m\left(  k\right)  -1},t_{0}\right)  =1-\varepsilon_{0}.
\]

\end{definition}

Notice that the previous definition does not depend on any t-norm. However,
when a t-norm of a specific class plays a role, then additional properties hold.

\begin{theorem}
\label{K38 - 24 th Non-Arch impies NC}Each non-Arquimedean KM-fuzzy metric
space $(X,M,\ast)$ whose t-norm $\ast$ is continuous at the $1$-boundary
satisfies the property $\mathcal{NC}$.
\end{theorem}

\begin{proof}
Suppose that $\{x_{n}\}\subseteq X$ is a sequence which is not $M$-Cauchy and
verifies%
\[
\lim_{n\rightarrow\infty}M\left(  x_{n},x_{n+1},t\right)  =1\quad\text{for all
}t>0.
\]
Then there are $\varepsilon_{0}\in\left(  0,1\right)  $, $t_{0}>0$ and two
partial subsequences $\{x_{n(k)}\}_{k\in\mathbb{N}}$ and $\{x_{m(k)}%
\}_{k\in\mathbb{N}}$ of $\{x_{n}\}$ such that $k<n\left(  k\right)  <m\left(
k\right)  <n\left(  k+1\right)  $ and%
\[
M\left(  x_{n\left(  k\right)  },x_{m\left(  k\right)  -1},t_{0}\right)
>1-\varepsilon_{0}\geq M\left(  x_{n\left(  k\right)  },x_{m\left(  k\right)
},t_{0}\right)  \quad\text{for all }k\in\mathbb{N}.
\]
As $(X,M,\ast)$ is a non-Archimedean KM-FMS, then, for all $k\in\mathbb{N}$,%
\begin{align}
1-\varepsilon_{0} &  \geq M\left(  x_{n\left(  k\right)  },x_{m\left(
k\right)  },t_{0}\right)  \geq M\left(  x_{n\left(  k\right)  },x_{m\left(
k\right)  -1},t_{0}\right)  \ast M\left(  x_{m\left(  k\right)  -1}%
,x_{m\left(  k\right)  },t_{0}\right)  \nonumber\\
&  \geq\left(  1-\varepsilon_{0}\right)  \ast M\left(  x_{m\left(  k\right)
-1},x_{m\left(  k\right)  },t_{0}\right)  .\label{K38 - 29 prop}%
\end{align}
Since $\lim_{k\rightarrow\infty}M\left(  x_{m\left(  k\right)  -1},x_{m\left(
k\right)  },t_{0}\right)  =1$ and $\ast$ is continuous at the $1$-boundary,
then%
\[
\lim_{k\rightarrow\infty}M\left(  x_{n\left(  k\right)  },x_{m\left(
k\right)  },t_{0}\right)  =1-\varepsilon_{0}.
\]
Taking into account that, by (\ref{K38 - 29 prop}),%
\[
\overset{a_{k}}{\overbrace{\,1-\varepsilon_{0}\,}}~\geq~\overset{c_{k}%
}{\overbrace{\,M\left(  x_{n\left(  k\right)  },x_{m\left(  k\right)
-1},t_{0}\right)  \,}}\ast\overset{d_{k}}{\overbrace{\,M\left(  x_{m\left(
k\right)  -1},x_{m\left(  k\right)  },t_{0}\right)  \,}}~\geq~\overset{e_{k}%
}{\overbrace{\,\left(  1-\varepsilon_{0}\right)  \ast M\left(  x_{m\left(
k\right)  -1},x_{m\left(  k\right)  },t_{0}\right)  \,}},
\]
for all $k\in\mathbb{N}$ and $\ast$ is continuous, Proposition
\ref{K38 - 18 propo cancellation} (applied with $b_{k}=1$ for all
$k\in\mathbb{N}$) guarantees that%
\[
\lim_{k\rightarrow\infty}M\left(  x_{n\left(  k\right)  },x_{m\left(
k\right)  -1},t_{0}\right)  =1-\varepsilon_{0},
\]
which means that%
\[
\lim_{k\rightarrow\infty}M\left(  x_{n\left(  k\right)  },x_{m\left(
k\right)  },t_{0}\right)  =\lim_{k\rightarrow\infty}M\left(  x_{n\left(
k\right)  },x_{m\left(  k\right)  -1},t_{0}\right)  =1-\varepsilon_{0}.
\]
Next, observe that, for all $k\in\mathbb{N}$,%
\begin{align*}
\overset{a_{k}}{\overbrace{\,M\left(  x_{n\left(  k\right)  },x_{m\left(
k\right)  },t_{0}\right)  \,}~} &  \geq~\overset{b_{k}}{\overbrace{\,M\left(
x_{n\left(  k\right)  },x_{n\left(  k\right)  -1},t_{0}\right)  \,}}%
\ast\,\overset{c_{k}}{\overbrace{\,M\left(  x_{n\left(  k\right)
-1},x_{m\left(  k\right)  -1},t_{0}\right)  \,}}\ast\overset{d_{k}}%
{\overbrace{\,M\left(  x_{m\left(  k\right)  -1},x_{m\left(  k\right)  }%
,t_{0}\right)  \,}}\\
&  \geq~M\left(  x_{n\left(  k\right)  },x_{n\left(  k\right)  -1}%
,t_{0}\right)  \ast M\left(  x_{n\left(  k\right)  -1},x_{n(k)},t_{0}\right)
\ast M\left(  x_{n\left(  k\right)  },x_{m\left(  k\right)  },t_{0}\right)  \\
&  \qquad\underset{e_{k}}{\underbrace{\qquad\qquad\ast M\left(  x_{m\left(
k\right)  },x_{m\left(  k\right)  -1},t_{0}\right)  \ast M\left(  x_{m\left(
k\right)  -1},x_{m\left(  k\right)  },t_{0}\right)  \qquad\qquad}}~.
\end{align*}
Clearly $\{a_{k}\}\rightarrow1-\varepsilon_{0}$, $\{b_{k}\}\rightarrow1$,
$\{d_{k}\}\rightarrow1$ and $\{e_{k}\}\rightarrow1-\varepsilon_{0}$. Thus,
Proposition \ref{K38 - 18 propo cancellation} guarantees again that%
\[
\lim_{k\rightarrow\infty}M\left(  x_{n\left(  k\right)  -1},x_{m\left(
k\right)  -1},t_{0}\right)  =1-\varepsilon_{0}.
\]
As a consequence,%
\[
\lim_{k\rightarrow\infty}M\left(  x_{n\left(  k\right)  },x_{m\left(
k\right)  },t_{0}\right)  =\lim_{k\rightarrow\infty}M\left(  x_{n\left(
k\right)  -1},x_{m\left(  k\right)  -1},t_{0}\right)  =1-\varepsilon_{0}.
\]
This completes the proof.
\end{proof}

\begin{corollary}
\label{K38 - 43 coro metric implies NC}If $(X,d)$ is a metric space, then
$(X,M^{d})$ satisfies the property $\mathcal{NC}$.
\end{corollary}

\begin{proof}
By item \ref{K38 - 44 propo Md non-Archimedean, item 2} of Proposition
\ref{K38 - 44 propo Md non-Archimedean, item 3}, $(X,M^{d},\ast_{p})$ is a
non-Archimedean KM-FMS and, as $\ast_{P}$ is continuos, Theorem
\ref{K38 - 24 th Non-Arch impies NC} ensures that $(X,M^{d})$ satisfies the
property $\mathcal{NC}$.
\end{proof}

\begin{example}
Let $(X,d)$ be a metric space for which there are $x_{0},y_{0},z_{0}\in X$
such that $d(x_{0},z_{0})>\max\{d(x_{0},y_{0}),d(y_{0},z_{0})\}$. By Corollary
\ref{K38 - 43 coro metric implies NC}, the fuzzy space $(X,M^{d})$ satisfies
the property $\mathcal{NC}$. Furthermore, item
\ref{K38 - 44 propo Md non-Archimedean, item 1} of Proposition
\ref{K38 - 44 propo Md non-Archimedean, item 3} ensures that $(X,M^{d}%
,\ast_{m})$ is a KM-FMS. However, taking into account the third item of the
same proposition, $(X,M^{d},\ast_{m})$ is not a non-Archimedean KM-FMS because
it does not satisfy the non-Archimedean property for $x_{0}$, $y_{0}$ and
$z_{0}$.
\end{example}

\section{\textbf{Fuzzy ample spectrum contractions}}

In this section we introduce two distinct notions of \emph{fuzzy ample
spectrum contraction} in the setting of KM-FMS. At a first sight, we can
believe that they are the natural performance in FMS of the properties that
define an ample spectrum contraction in metric spaces. However, FMS are
distinct in nature to metric spaces: for instance, they are much varied than
metric spaces. As we advise in the introduction, we will work on a KM-FMS.
Their main advantage is that they are more general than GV-FMS, but they also
have a great drawback: the statements of the fixed point theorems need more
technical hypotheses and their corresponding proofs are more difficult. The
reader can easily deduce how the proofs would be easier in a GV-FMS.
Furthermore, rather than working in a non-Archimedean FMS, we will employ
KM-FMS that only satisfy the property $\mathcal{NC}$ and whose t-norms are
continuous at the $1$-boundary (which are more general).

We give the following definitions in the context of fuzzy spaces. Throughout
this section, let $\left(  X,M\right)  $ be a fuzzy space, let $T:X\rightarrow
X$ be a self-mapping, let $\mathcal{S}$ be a binary relation on $X$, let
$B\subseteq\mathbb{R}^{2}$ be a subset and let $\theta:B\rightarrow\mathbb{R}$
be a function. Let%
\[
\mathcal{M}_{T}=\left\{  \,(M\left(  Tx,Ty,t\right)  ,M\left(  x,y,t\right)
)\in\mathbb{I}\times\mathbb{I}:x,y\in X,~x\mathcal{S}^{\ast}y,~Tx\mathcal{S}%
^{\ast}Ty,~t\in\left[  0,\infty\right)  \,\right\}  .
\]

\subsection{(\textbf{Type-1) Fuzzy ample spectrum contractions}}

The notion of fuzzy ample spectrum contraction directly depend on a very
singular kind of sequences of pairs of distance distribution functions.

\begin{definition}
\label{K38 - 55 def T,S,M sequence}Let $\{\phi_{n}\}$ and $\{\psi_{n}\}$ be
two sequences of functions $\phi_{n},\psi_{n}:\left[  0,\infty\right)
\rightarrow\mathbb{I}$ . We say that $\{\left(  \phi_{n},\psi_{n}\right)  \}$
is a {$(T,\mathcal{S}^{\ast},M)$\emph{-sequence}} if there exist two sequences
$\{x_{n}\},\{y_{n}\}\subseteq X$ such that%
\[
x_{n}\mathcal{S}^{\ast}y_{n},\quad Tx_{n}\mathcal{S}^{\ast}Ty_{n},\quad
\phi_{n}\left(  t\right)  =M(Tx_{n},Ty_{n},t)\quad\text{and}\quad\psi
_{n}\left(  t\right)  =M(x_{n},y_{n},t)
\]
for all $n\in\mathbb{N}$ and all $t>0$.
\end{definition}

\begin{definition}
\label{K38 - 54 def FASC}A mapping $T:X\rightarrow X$ is said to be a
\emph{fuzzy ample spectrum contraction w.r.t. }$\left(  M,\mathcal{S}%
,\theta\right)  $ if the following four conditions are fulfilled.

\begin{description}
\item[$(\mathcal{F}_{1})$] $B$ is nonempty and $\mathcal{M}_{T}\subseteq B$.

\item[$(\mathcal{F}_{2})$] If $\{x_{n}\}\subseteq X$ is a Picard $\mathcal{S}
$-strictly-increasing sequence of $T$ such that
\[
\theta\left(  M\left(  x_{n+1},x_{n+2},t\right)  ,M\left(  x_{n}%
,x_{n+1},t\right)  \right)  \geq0\quad\text{for all }n\in\mathbb{N}\text{ and
all }t>0,
\]
then $\lim_{n\rightarrow\infty}M\left(  x_{n},x_{n+1},t\right)  =1$ for all
$t>0$.

\item[$(\mathcal{F}_{3})$] If $\{\left(  \phi_{n},\psi_{n}\right)  \}$ is a
$(T,\mathcal{S}^{\ast},M)$-sequence and $t_{0}>0$ are such that $\{\phi
_{n}\left(  t_{0}\right)  \}$ and $\{\psi_{n}\left(  t_{0}\right)  \}$
converge to the same limit $L\in\mathbb{I}$ and verifying that $L>\phi
_{n}\left(  t_{0}\right)  $ and $\theta(\phi_{n}\left(  t_{0}\right)
,\psi_{n}\left(  t_{0}\right)  )\geq0$ for all $n\in\mathbb{N}_{0}$, then
$L=1$.

\item[$(\mathcal{F}_{4})$] ${\theta}${$\left(  M(Tx,Ty,t),M(x,y,t)\right)
\geq0\quad$for all }$t>0$ and all {\ $x,y\in X\quad$such that $x\mathcal{S}%
^{\ast}y$ and $Tx\mathcal{S}^{\ast}Ty$.}
\end{description}
\end{definition}

In some cases, we will also consider the following properties.

\begin{description}
\item[$(\mathcal{F}_{2}^{\prime})$] If $\{x_{n}\},\{y_{n}\}\subseteq X$ are
two $T$-Picard sequences such that%
\[
x_{n}\mathcal{S}^{\ast}y_{n}\quad\text{and}\quad{\theta(M\left(
x_{n+1},y_{n+1},t\right)  ,M\left(  x_{n},y_{n},t\right)  )\geq0}%
\quad\text{for all }n\in\mathbb{N}\text{ and all }t>0,
\]
then $\lim_{n\rightarrow\infty}{M}\left(  x_{n},y_{n},t\right)  =1$ for all
$t>0$.

\item[$(\mathcal{F}_{5})$] If $\{\left(  \phi_{n},\psi_{n}\right)  \}$ is a
$(T,\mathcal{S}^{\ast},M)$-sequence such that $\{\psi_{n}\left(  t\right)
\}\rightarrow1$ for all $t>0$ and $\theta(\phi_{n}(t),\psi_{n}(t))\geq0$ for
all $n\in\mathbb{N}$ and all $t>0$, then $\{\phi_{n}\left(  t\right)
\}\rightarrow1$ for all $t>0$.
\end{description}

Many of the remarks that were given in the context of metric spaces for ample
spectrum contractions can now be repeated. In particular, we highlight the
following ones.

\begin{remark}
\label{K38 - 30 rem def FAEC}

\begin{enumerate}
\item \label{K38 - 30 rem def FAEC, item 1}Although the set $B$ in which is
defined the function $\theta:B\rightarrow\mathbb{R}$ can be greater than
$\mathbb{I}\times\mathbb{I}$, for our purposes, we will only be interested in
the values of $\theta$ when its arguments belong to $\mathbb{I}\times
\mathbb{I}$. Hence it is sufficient to assume that $B\subseteq\mathbb{I}%
\times\mathbb{I}$.

\item \label{K38 - 30 rem def FAEC, item 2}If the function $\theta$ satisfies
$\theta\left(  t,s\right)  \leq t-s$ for all $\left(  t,s\right)  \in
B\cap(\mathbb{I}\times\mathbb{I})$, then property $(\mathcal{F}_{5})$ holds.
It follows from the fact that%
\[
0\leq\theta(\phi_{n}(t),\psi_{n}(t))\leq\phi_{n}(t)-\psi_{n}(t)\quad
\Rightarrow\quad\psi_{n}(t)\leq\phi_{n}(t)\leq1.
\]

\item By choosing $y_{n}=x_{n+1}$ for all $n\in\mathbb{N}$, it can be proved
that $(\mathcal{F}_{2}^{\prime})$ implies $(\mathcal{F}_{2})$.
\end{enumerate}
\end{remark}

In the following result we check that a great family of ample spectrum
contractions are, in fact, fuzzy ample spectrum contractions.

\begin{theorem}
\label{K38 - 56 th ASC implies FASC}Let $(X,d)$ be a metric space and let
$T:X\rightarrow X$ be an ample spectrum contraction w.r.t. $\mathcal{S}$ and
$\varrho$, where $\varrho:\left[  0,\infty\right)  \times\left[
0,\infty\right)  \rightarrow\mathbb{R}$. Suppose that the function $\varrho$
satisfies the following property:%
\begin{equation}
t,s>0,\quad\varrho\left(  t,s\right)  \geq0\quad\Rightarrow\quad\left[
~\varrho\left(  \frac{t}{r},\frac{s}{r}\right)  \geq0~~\text{for all
}r>0~\right]  . \label{K38 - 51 prop}%
\end{equation}
Let define:%
\[
\theta_{\varrho}:\left(  0,1\right]  \times\left(  0,1\right]  \rightarrow
\mathbb{R},\qquad\theta_{\varrho}\left(  t,s\right)  =\varrho\left(
\frac{1-t}{t},\frac{1-s}{s}\right)  \quad\text{for all }t,s\in\left(
0,1\right]  .
\]
Then $T$ is a fuzzy ample spectrum contraction w.r.t. $(M^{d},\mathcal{S}%
,\theta_{\varrho})$. Furthermore, if property $(\mathcal{B}_{5})$
(respectively, $(\mathcal{B}_{2}^{\prime})$) holds, then property
$(\mathcal{F}_{5})$ (respectively, $(\mathcal{F}_{2}^{\prime})$) also holds.
\end{theorem}

\begin{proof}
First of all, we observe that, for all $x,y\in X$ and all $t>0$,%
\begin{align}
&  \theta_{\varrho}(M^{d}(Tx,Ty,t),M^{d}(x,y,t))=\theta_{\varrho}\left(
\frac{t}{t+d(Tx,Ty)},\frac{t}{t+d(x,y)}\right) \nonumber\\
&  \qquad=\varrho\left(  \frac{\,1-\dfrac{t}{t+d(Tx,Ty)}\,}{\dfrac
{t}{t+d(Tx,Ty)}},\frac{\,1-\dfrac{t}{t+d(x,y)}\,}{\dfrac{t}{t+d(x,y)}}\right)
=\varrho\left(  \frac{d(Tx,Ty)}{t},\frac{d(x,y)}{t}\right)  .
\label{K38 - 50 prop}%
\end{align}
Next we check all properties that define a fuzzy ample spectrum contraction.

$(\mathcal{F}_{1})$ $B=\left(  0,1\right]  \times\left(  0,1\right]  $ is
nonempty and $\mathcal{M}_{T}\subseteq B$.

$(\mathcal{F}_{2})$ Let $\{x_{n}\}\subseteq X$ be a Picard $\mathcal{S}%
$-strictly-increasing sequence of $T$ such that
\[
\theta_{\varrho}(M^{d}\left(  x_{n+1},x_{n+2},t\right)  ,M^{d}\left(
x_{n},x_{n+1},t\right)  )\geq0\quad\text{for all }n\in\mathbb{N}\text{ and all
}t>0.
\]
In particular, for $t=1$, using (\ref{K38 - 50 prop}), for each $n\in
\mathbb{N}$,%
\[
0\leq\theta_{\varrho}(M^{d}\left(  x_{n+1},x_{n+2},1\right)  ,M^{d}\left(
x_{n},x_{n+1},1\right)  )=\varrho(d(x_{n+1},x_{n+2}),d(x_{n},x_{n+1})).
\]
As $T$ is an ample spectrum contraction w.r.t. $\mathcal{S}$ and $\varrho$,
axiom $(\mathcal{B}_{2})$ implies that $\{d(x_{n},x_{n+1})\}\rightarrow0$, so
$\lim_{n\rightarrow\infty}M^{d}\left(  x_{n},x_{n+1},t\right)  =1$ for all
$t>0$.

$(\mathcal{F}_{3})$~Let $\{\left(  \phi_{n},\psi_{n}\right)  \}$ be a
$(T,\mathcal{S}^{\ast},M^{d})$-sequence and let $t_{0}>0$ be such that
$\{\phi_{n}\left(  t_{0}\right)  \}$ and $\{\psi_{n}\left(  t_{0}\right)  \}$
converge to the same limit $L\in\mathbb{I}$ and verifying that $L>\phi
_{n}\left(  t_{0}\right)  $ and $\theta_{\varrho}(\phi_{n}\left(
t_{0}\right)  ,\psi_{n}\left(  t_{0}\right)  )\geq0$ for all $n\in
\mathbb{N}_{0}$. Let $\{x_{n}\},\{y_{n}\}\subseteq X$ be some sequences such
that, for all $n\in\mathbb{N}$ and all $t>0$,%
\[
x_{n}\mathcal{S}^{\ast}y_{n},\quad Tx_{n}\mathcal{S}^{\ast}Ty_{n},\quad
\phi_{n}\left(  t\right)  =M^{d}(Tx_{n},Ty_{n},t)\quad\text{and}\quad\psi
_{n}\left(  t\right)  =M^{d}(x_{n},y_{n},t).
\]
Let define $a_{n}=d(Tx_{n},Ty_{n})$ and $b_{n}=d(x_{n},y_{n})$ for all
$n\in\mathbb{N}$. Then $\{\left(  a_{n},b_{n}\right)  \}\subseteq\left[
0,\infty\right)  \times\left[  0,\infty\right)  $ is a $(T,\mathcal{S}^{\ast
})$-sequence. Furthermore, for all $n\in\mathbb{N}$,%
\begin{align*}
0  &  \leq\theta_{\varrho}(\phi_{n}\left(  t_{0}\right)  ,\psi_{n}\left(
t_{0}\right)  )=\theta_{\varrho}(M^{d}(Tx_{n},Ty_{n},t_{0}),M^{d}(x_{n}%
,y_{n},t_{0}))\\[0.2cm]
&  =\varrho\left(  \frac{d(Tx_{n},Ty_{n})}{t_{0}},\frac{d(x_{n},y_{n})}{t_{0}%
}\right)  =\varrho\left(  \frac{a_{n}}{t_{0}},\frac{b_{n}}{t_{0}}\right)  .
\end{align*}
By (\ref{K38 - 51 prop}), using $r=1/t_{0}$, we deduce that, for all
$n\in\mathbb{N}$,%
\[
0\leq\varrho\left(  \frac{\,\frac{a_{n}}{t_{0}}\,}{\frac{1}{t_{0}}}%
,\frac{\,\frac{b_{n}}{t_{0}}\,}{\frac{1}{t_{0}}}\right)  =\varrho\left(
a_{n},b_{n}\right)  .
\]
Since $L>\phi_{n}\left(  t_{0}\right)  $, then $L>0$. Let $\varepsilon
_{0}=t_{0}\frac{1-L}{L}\geq0$. Thus, as $\{\phi_{n}\left(  t_{0}\right)
\}\rightarrow L$ and $\{\psi_{n}\left(  t_{0}\right)  \}\rightarrow L$, then%
\[
\left\{  \frac{t_{0}}{t_{0}+d(Tx_{n},Ty_{n})}\right\}  _{n\in\mathbb{N}%
}\rightarrow L\quad\text{and}\quad\left\{  \frac{t_{0}}{t_{0}+d(x_{n},y_{n}%
)}\right\}  _{n\in\mathbb{N}}\rightarrow L.
\]
This is equivalent to%
\[
\{a_{n}\}=\left\{  d(Tx_{n},Ty_{n})\right\}  _{n\in\mathbb{N}}\rightarrow
t_{0}\frac{1-L}{L}=\varepsilon_{0}\quad\text{and}\quad\{b_{n}\}=\left\{
d(x_{n},y_{n})\right\}  _{n\in\mathbb{N}}\rightarrow t_{0}\frac{1-L}%
{L}=\varepsilon_{0}.
\]
Notice that%
\begin{align*}
\frac{t_{0}}{t_{0}+d(Tx_{n},Ty_{n})}=\phi_{n}\left(  t_{0}\right)  <L\quad &
\Leftrightarrow\quad\frac{1}{t_{0}+a_{n}}<\frac{L}{t_{0}}\quad\Leftrightarrow
\quad\frac{t_{0}}{L}<t_{0}+a_{n}\\[0.2cm]
&  \Leftrightarrow\quad\frac{t_{0}}{L}-t_{0}<a_{n}\quad\Leftrightarrow
\quad\varepsilon_{0}=t_{0}\frac{1-L}{L}<a_{n}.
\end{align*}
Taking into account that $T$ is an ample spectrum contraction w.r.t.
$\mathcal{S}$ and $\varrho$, condition $(\mathcal{B}_{3})$ ensures that
$\varepsilon_{0}=0$, which leads to $L=1$.

$(\mathcal{F}_{4})$ Let{\ }$t>0$ and let {$x,y\in X$ be such that
$x\mathcal{S}^{\ast}y$ and $Tx\mathcal{S}^{\ast}Ty$. Therefore, }$d(x,y)>0$
and $d(Tx,Ty)>0$.{\ As} $T$ is an ample spectrum contraction w.r.t.
$\mathcal{S}$ and $\varrho$, property $(\mathcal{B}_{4})$ guarantees that
$\varrho\left(  d(Tx,Ty),d(x,y)\right)  \geq0$. By (\ref{K38 - 51 prop}), it
follows that%
\[
\varrho\left(  \frac{d(Tx,Ty)}{r},\frac{d(x,y)}{r}\right)  \geq0~~\text{for
all }r>0.
\]
As a consequence, by (\ref{K38 - 50 prop}),%
\[
\theta_{\varrho}({M^{d}(Tx,Ty,t),M^{d}(x,y,t))=\varrho\left(  \frac
{d(Tx,Ty)}{t},\frac{d(x,y)}{t}\right)  \geq0.}%
\]

$(\mathcal{B}_{5})\Rightarrow(\mathcal{F}_{5})$. Suppose that the property
$(\mathcal{B}_{5})$ holds. Let $\{\left(  \phi_{n},\psi_{n}\right)  \}$ be a
$(T,\mathcal{S}^{\ast},M)$-sequence such that $\{\psi_{n}\left(  t\right)
\}\rightarrow1$ for all $t>0$ and $\theta_{\varrho}(\phi_{n}(t),\psi
_{n}(t))\geq0$ for all $n\in\mathbb{N}$ and all $t>0$. Let $\{x_{n}%
\},\{y_{n}\}\subseteq X$ be some sequences such that, for all $n\in\mathbb{N}$
and all $t>0$,%
\[
x_{n}\mathcal{S}^{\ast}y_{n},\quad Tx_{n}\mathcal{S}^{\ast}Ty_{n},\quad
\phi_{n}\left(  t\right)  =M^{d}(Tx_{n},Ty_{n},t)\quad\text{and}\quad\psi
_{n}\left(  t\right)  =M^{d}(x_{n},y_{n},t).
\]
Using $t=1$, for all $n\in\mathbb{N}$,%
\[
0\leq\theta_{\varrho}(\phi_{n}(1),\psi_{n}(1))=\theta_{\varrho}(M^{d}%
(Tx_{n},Ty_{n},1),M^{d}(x_{n},y_{n},1))=\varrho\left(  d(Tx_{n},Ty_{n}%
),d(x_{n},y_{n})\right)  .
\]
Notice that $\{M^{d}(x_{n},y_{n},t)\}=\{\psi_{n}\left(  t\right)
\}\rightarrow1$ for all $t>0$ is equivalent to say that $\{d(x_{n}%
,y_{n})\}\rightarrow0$. Taking into account that we suppose that property
$(\mathcal{B}_{5})$ holds, then $\{d(Tx_{n},Ty_{n})\}\rightarrow0$, so
$\{M^{d}(Tx_{n},Ty_{n},t)\}=\{\phi_{n}\left(  t\right)  \}\rightarrow1$ for
all $t>0$.

$(\mathcal{B}_{2}^{\prime})\Rightarrow(\mathcal{F}_{2}^{\prime})$. Let
$\{x_{n}\},\{y_{n}\}\subseteq X$ be two $T$-Picard sequences such that%
\[
x_{n}\mathcal{S}^{\ast}y_{n}\quad\text{and}\quad{\theta}_{\varrho}{(M}%
^{d}{\left(  x_{n+1},y_{n+1},t\right)  ,M}^{d}{\left(  x_{n},y_{n},t\right)
)\geq0}\quad\text{for all }n\in\mathbb{N}\text{ and all }t>0.
\]
Using $t=1$, for all $n\in\mathbb{N}$,%
\[
0\leq{\theta}_{\varrho}{(M}^{d}{\left(  x_{n+1},y_{n+1},1\right)  ,M}%
^{d}{\left(  x_{n},y_{n},1\right)  )}=\varrho\left(  d(Tx_{n},Ty_{n}%
),d(x_{n},y_{n})\right)  .
\]
As we suppose that property $(\mathcal{B}_{2}^{\prime})$ holds, then
$\{d(x_{n},y_{n})\}\rightarrow0$ which means that $\{{M}^{d}{\left(
x_{n},y_{n},t\right)  }\}\rightarrow1$ for all $t>0$.
\end{proof}

\begin{example}
If $T$ is a Banach contraction, then there is $\lambda\in\left[  0,1\right)  $
such that $d(Tx,Ty)\leq\lambda d(x,y)$ for all $x,y\in X$. In this case, $T$
is an ample spectrum contraction associated to the function $\varrho_{\lambda
}(t,s)=\lambda s-t$ for all $t,s\geq0$. In such a case, for all $t,s,r>0$,%
\[
\varrho_{\lambda}\left(  \frac{t}{r},\frac{s}{r}\right)  =\lambda\frac{s}%
{r}-\frac{t}{r}=\frac{\lambda s-t}{r}=\frac{\varrho_{\lambda}(t,s)}{r},
\]
which means that property (\ref{K38 - 51 prop}) holds.
\end{example}

In this general framework we introduce one of our main results. Notice that in
the following result we assume that the t-norm $\ast$ is continuous at the
$1$-boundary, but it has not to be continuous on the whole space
$\mathbb{I}\times\mathbb{I}$.

\begin{theorem}
\label{K38 - 31 th main FAEC}Let $\left(  X,M,\ast\right)  $ be a KM-FMS
endowed with a transitive binary relation $\mathcal{S}$ and let
$T:X\rightarrow X$ be an $\mathcal{S}$-nondecreasing fuzzy ample spectrum
contraction with respect to $\left(  M,\mathcal{S},\theta\right)  $. Suppose
that $\ast$ is continuous at the $1$-boundary, $T(X)$ is
$(\mathcal{S},M)$-strictly-increasing-precomplete and there exists a point $x_{0}\in X$
such that $x_{0}\mathcal{S}Tx_{0}$. Also assume that, at least, one of the
following conditions is fulfilled:

\begin{description}
\item[$\left(  a\right)  $] $T$ is $\mathcal{S}$-strictly-increasing-continuous.

\item[$\left(  b\right)  $] $\left(  X,M\right)  $ is $\mathcal{S}%
$-strictly-increasing-regular and condition $(\mathcal{F}_{5})$ holds.

\item[$\left(  c\right)  $] $\left(  X,M\right)  $ is $\mathcal{S}%
$-strictly-increasing-regular and $\theta\left(  t,s\right)  \leq t-s$ for all
$\left(  t,s\right)  \in B\cap(\mathbb{I}\times\mathbb{I})$.
\end{description}

If $\left(  X,M\right)  $\ satisfies the property $\mathcal{NC}$, then the
Picard sequence of $T~$based on $x_{0}$ converges to a fixed point of $T$. In
particular, $T$ has at least a fixed point.
\end{theorem}

\begin{proof}
Let $x_{0}\in X$ be the point such that $x_{0}\mathcal{S}Tx_{0}$ and let
$\{x_{n}\}$ be the Picard sequence of $T$ starting from $x_{0}$. If there is
some $n_{0}\in\mathbb{N}_{0}$ such that $x_{n_{0}}=x_{n_{0}+1}$, then
$x_{n_{0}}$ is a fixed point of $T$, and the proof is finished. On the
contrary, suppose that $x_{n}\neq x_{n+1}$ for all $n\in\mathbb{N}_{0}$. As
$x_{0}\mathcal{S}Tx_{0}=x_{1}$ and $T$ is $\mathcal{S}$-nondecreasing, then
$x_{n}\mathcal{S}x_{n+1}$ for all $n\in\mathbb{N}_{0}$. In fact, as
$\mathcal{S}$ is transitive, then%
\[
x_{n}\mathcal{S}x_{m}\quad\text{for all }n,m\in\mathbb{N}_{0}\text{ such that
}n<m,
\]
which means that $\{x_{n}\}$ is an $\mathcal{S}$-nondecreasing sequence. Since
$T$ is a fuzzy ample spectrum contraction with respect to $\left(
M,\mathcal{S},\theta\right)  $ and $x_{n}\mathcal{S}^{\ast}x_{n+1}$ and
$Tx_{n}=x_{n+1}\mathcal{S}^{\ast}x_{n+2}=Tx_{n+1}$, then%
\[
{\theta\left(  M(x_{n+1},x_{n+2},t),M(x_{n},x_{n+1},t)\right)  ={\theta\left(
M(Tx_{n},Tx_{n+1},t),M(x_{n},x_{n+1},t)\right)  }\geq0}%
\]
{for all }$t>0$ and all $n\in\mathbb{N}_{0}${. Axiom }$(\mathcal{F}_{2})$
guarantees that%
\begin{equation}
\lim_{n\rightarrow\infty}M\left(  x_{n},x_{n+1},t\right)  \rightarrow
1\quad\text{for all }t>0. \label{K38 - 17 prop}%
\end{equation}
By Proposition \ref{K38 - 22 lem infinite or almost periodic}, the $T$-Picard
sequence $\{x_{n}\}$ is infinite or almost periodic. If $\{x_{n}\}$ is almost
periodic, then there are $n_{0},m_{0}\in\mathbb{N}$ such that $n_{0}<m_{0}$
and $x_{n_{0}}=x_{m_{0}}$. In such a case, Proposition
\ref{K40 21 propo either infinite or almost-constant fuzzy} guarantees that
there is $\ell_{0}\in\mathbb{N}$ and $z\in X$ such that $x_{n}=z$ for all
$n\geq\ell_{0}$, so $z$ is a fixed point of $T$, and the proof is finished. On
the contrary case, suppose that $\{x_{n}\}$ is infinite, that is, $x_{n}\neq
x_{m}$ for all $n\neq m$. This means that $\{x_{n}\}$ is an $\mathcal{S}%
$-strictly-increasing sequence because $x_{n}\mathcal{S}x_{m}$ and $x_{n}\neq
x_{m}$, that is, $x_{n}\mathcal{S}^{\ast}x_{m}$ for all $n<m$.

Next we prove that $\{x_{n}\}$ is an $M$-Cauchy sequence by contradiction.
Since we suppose that $\left(  X,M\right)  $ satisfies the property
$\mathcal{NC}$ and $\{x_{n}\}$ is a sequence which is not $M$-Cauchy but it
satisfies (\ref{K38 - 17 prop}), there are $\varepsilon_{0}\in\left(
0,1\right)  $ and $t_{0}>0$ and two partial subsequences $\{x_{n(k)}%
\}_{k\in\mathbb{N}}$ and $\{x_{m(k)}\}_{k\in\mathbb{N}}$ of $\{x_{n}\}$ such
that, for all $k\in\mathbb{N}$, $k<n\left(  k\right)  <m\left(  k\right)
<n\left(  k+1\right)  $,%
\[
M\left(  x_{n\left(  k\right)  },x_{m\left(  k\right)  -1},t_{0}\right)
\geq1-\varepsilon_{0}>M\left(  x_{n\left(  k\right)  },x_{m\left(  k\right)
},t_{0}\right)
\]
and also%
\begin{equation}
\lim_{k\rightarrow\infty}M\left(  x_{n\left(  k\right)  },x_{m\left(
k\right)  },t_{0}\right)  =\lim_{k\rightarrow\infty}M\left(  x_{n\left(
k\right)  -1},x_{m\left(  k\right)  -1},t_{0}\right)  =1-\varepsilon_{0}.
\label{K38 - 21 prop}%
\end{equation}
Let define $L=1-\varepsilon_{0}\in\left(  0,1\right)  $ and, for all $t>0$ and
all $k\in\mathbb{N}$,%
\[
\phi_{k}\left(  t\right)  =M(x_{n\left(  k\right)  },x_{m\left(  k\right)
},t)\quad\text{and}\quad\psi_{k}\left(  t\right)  =M(x_{n\left(  k\right)
-1},x_{m\left(  k\right)  -1},t).
\]
We claim that $\{(\phi_{k},\psi_{k})\}_{k\in\mathbb{N}}$, $t_{0}$ and $L$
satisfy all hypotheses in condition $(\mathcal{F}_{3})$. On the one hand,
$\{(\phi_{k},\psi_{k})\}$ is a {$(T,\mathcal{S}^{\ast},M)$}-sequence because
$\phi_{k}=M(x_{n\left(  k\right)  },x_{m\left(  k\right)  },\cdot
)=M(Tx_{n\left(  k\right)  -1},Tx_{m\left(  k\right)  -1},\cdot)$, $\psi
_{k}=M(x_{n\left(  k\right)  -1},x_{m\left(  k\right)  -1},\cdot)$,
$x_{n\left(  k\right)  -1}\mathcal{S}^{\ast}x_{m\left(  k\right)  -1}$ and
$Tx_{n\left(  k\right)  -1}=x_{n\left(  k\right)  }\mathcal{S}^{\ast
}x_{m\left(  k\right)  }=Tx_{m\left(  k\right)  -1}$ for all $k\in\mathbb{N}$.
Furthermore, $L=1-\varepsilon_{0}>M(x_{n\left(  k\right)  },x_{m\left(
k\right)  },t_{0})=\phi_{k}\left(  t_{0}\right)  $ for all $k\in\mathbb{N}$.
Also (\ref{K38 - 21 prop}) means that%
\[
\lim_{k\rightarrow\infty}\phi_{k}\left(  t_{0}\right)  =\lim_{k\rightarrow
\infty}\psi_{k}\left(  t_{0}\right)  =1-\varepsilon_{0}=L.
\]
As $T$ is a fuzzy ample spectrum contraction with respect to $\left(
M,\mathcal{S},\theta\right)  $, condition $(\mathcal{F}_{3})$ guarantees that
$1-\varepsilon_{0}=L=1$, which is a contradiction because $\varepsilon_{0}>0
$. This contradiction shows that $\{x_{n}\}$ is an $M$-Cauchy sequence in
$(X,M)$. Since $\{x_{n+1}=Tx_{n}\}\subseteq TX$ and $T(X)$ is $(\mathcal{S},M)
$-strictly-increasing-precomplete, then there exists a set $Z$ such that
$TX\subseteq Z\subseteq X$ and $Z$ is $(\mathcal{S},M)$%
-strictly-increasing-complete. Since $\{x_{n+1}=Tx_{n}\}\subseteq TX\subseteq
Z$ and $\{x_{n}\}$ is an $\mathcal{S}$-strictly-increasing $M$-Cauchy
sequence, then there is $z\in Z\subseteq X$ such that $\{x_{n}\}$
$M$-converges to $z$. It only remain to prove that, under any of the
conditions $(a)$, $(b)$ or $(c)$, $z$ is a fixed point of $T$.

\begin{description}
\item[$(a)$] Suppose that $T$ is $\mathcal{S}$-strictly-increasing-continuous.
As $\{x_{n}\}$ is $\mathcal{S}$-strictly-increasing and $\{x_{n}\}$
$M$-converges to $z$, then $\{Tx_{n}\}$ $M$-converges to $Tz$. However, as
$Tx_{n}=x_{n+1}$ for all $n\in\mathbb{N}$ and the $M$-limit in a KM-FMS is
unique, then $Tz=z$, that is, $z$ is a fixed point of $T$.

\item[$(b)$] Suppose that $\left(  X,M\right)  $ is $\mathcal{S}%
$-strictly-increasing-regular and condition $(\mathcal{F}_{5})$ holds. In this
case, since $\{x_{n}\}$ is an $\mathcal{S}$-strictly-increasing sequence such
that $\{x_{n}\}\rightarrow z\in X$, it follows that $x_{n}\mathcal{S}z$ for
all $n\in\mathbb{N}$. Taking into account that the sequence $\{x_{n}\}$ is
infinite, then there is $n_{0}\in\mathbb{N}$ such that $x_{n}\neq z$ and
$x_{n}\neq Tz$ for all $n\geq n_{0}-1$. Moreover, as $T$ is $\mathcal{S}%
$-nondecreasing, then $x_{n+1}=Tx_{n}\mathcal{S}Tz$, which means that
$x_{n}\mathcal{S}^{\ast}z$ and $x_{n}\mathcal{S}^{\ast}Tz$ for all $n\geq
n_{0}$. Using that $T$ is a fuzzy ample spectrum contraction with respect to
$\left(  M,\mathcal{S},\theta\right)  $, condition $(\mathcal{F}_{4})$ implies
that%
\[
\theta\left(  M(x_{n+1},Tz,t),M(x_{n},z,t)\right)  =\theta\left(
M(Tx_{n},Tz,t),M(x_{n},z,t)\right)  \geq0
\]
for all $n\geq n_{0}$ and all $t>0$. Taking into account that $\{M(x_{n}%
,z,t)\}_{n\geq n_{0}}\rightarrow1$ for all $t>0$, assumption $(\mathcal{F}%
_{5})$ applied to the $(T,\mathcal{S}^{\ast},M)$-sequence%
\[
\left\{  \,\left(  \,\phi_{n}=M(Tx_{n},Tz,\cdot),\,\psi_{n}=M(x_{n}%
,z,\cdot)\,\right)  \,\right\}  _{n\geq n_{0}}%
\]
leads to $\{M(x_{n+1},Tz,t)\}_{n\geq n_{0}}\rightarrow1$ for all $t>0$, that
is, $\{x_{n}\}_{n\geq n_{0}}$ $M$-converges to $Tz$. As the $M$-limit in a
KM-FMS is unique, then $Tz=z$, so $z$ is a fixed point of $T$.

\item[$(c)$] Suppose that $\left(  X,M\right)  $ is $\mathcal{S}%
$-strictly-increasing-regular and $\theta\left(  t,s\right)  \leq t-s$ for all
$\left(  t,s\right)  \in B\cap(\mathbb{I}\times\mathbb{I})$. This case follows
from $(b)$ taking into account item \ref{K38 - 30 rem def FAEC, item 2} of
Remark \ref{K38 - 30 rem def FAEC}.
\end{description}
\end{proof}

Theorem \ref{K38 - 31 th main FAEC} guarantees the existence of fixed points
of $T$. In the following result we describe some additional assumptions in
order to ensure that such fixed point is unique.

\begin{theorem}
\label{K38 - 31 th main FAEC uniqueness}Under the hypotheses of Theorem
\ref{K38 - 31 th main FAEC}, assume that the property $(\mathcal{F}%
_{2}^{\prime}) $ is fulfilled and that each pair of fixed points
$x,y\in\operatorname*{Fix}(T)$ of $T$ is associated to another point $z\in X$
which is, at the same time, $\mathcal{S}$-comparable to $x$ and to $y$. Then
$T$ has a unique fixed point.
\end{theorem}

\begin{proof}
Let $x,y\in\operatorname*{Fix}(T)$ be two fixed points of $T$. By hypothesis,
there exists $z_{0}\in X$ such that $z_{0}$ is, at the same time,
$\mathcal{S}$-comparable to $x$ and $\mathcal{S}$-comparable to $y$. We claim
that the $T$-Picard sequence $\{z_{n}\}$ of $T$ starting from $z_{0}$
$M$-converges, at the same time, to $x$ and to $y$ (so we will deduce that
$x=y$). We check the first statement. To prove it, we consider two possibilities.

\begin{itemize}
\item Suppose that there is $n_{0}\in\mathbb{N}$ such that $z_{n_{0}}=x$. In
this case, $z_{n_{0}+1}=Tz_{n_{0}}=Tx=x$. Repeating this argument, $z_{n}=x$
for all $n\geq n_{0}$, so $\{z_{n}\}$ $M$-converges to $x$.

\item Suppose that $z_{n}\neq x$ for all $n\in\mathbb{N}$. Since $z_{0}$ is
$\mathcal{S}$-comparable to $x$, assume, for instance, that $z_{0}\mathcal{S}x
$ (the case $x\mathcal{S}z_{0}$ is similar). As $z_{0}\mathcal{S}x$, $T$ is
$\mathcal{S}$-nondecreasing and $Tx=x$, then $z_{n}\mathcal{S}x$ for all
$n\in\mathbb{N}$. Therefore $z_{n}\mathcal{S}^{\ast}x$ and $Tz_{n}%
\mathcal{S}^{\ast}Tx$ for all $n\in\mathbb{N}$. Using the contractivity
condition $(\mathcal{F}_{4})$, for all $n\in\mathbb{N}$ and all $t>0$,%
\[
0\leq\theta(M(Tz_{n},Tx,t),M(z_{n},x,t))=\theta(M(T^{n+1}z_{0},T^{n+1}%
x,t),M(T^{n}z_{0},T^{n}x,t)).
\]
It follows from $(\mathcal{F}_{2}^{\prime})$ that $\{M(z_{n},x,t)\}=\{M(T^{n}%
z_{0},T^{n}x,t)\}\rightarrow1$ for all $t>0$, that is, $\{z_{n}\}$
$M$-converges to $x$.
\end{itemize}

In any case, $\{z_{n}\}\rightarrow x$ and, similarly, $\{z_{n}\}\rightarrow y
$, so $x=y$ and $T$ has a unique fixed point.
\end{proof}

\begin{corollary}
Under the hypotheses of Theorem \ref{K38 - 31 th main FAEC}, assume that
condition $(\mathcal{F}_{2}^{\prime})$ holds and each two fixed points of $T$
are $\mathcal{S}$-comparable. Then $T$ has a unique fixed point.
\end{corollary}

Immediately we can deduce that Theorems \ref{K38 - 31 th main FAEC} and
\ref{K38 - 31 th main FAEC uniqueness} remains true if we replace the
hypothesis that $(X,M)$ satisfies the property $\mathcal{NC}$ by the fact that
$(X,M,\ast)$ is a non-Archimedean KM-FMS (recall Theorem
\ref{K38 - 24 th Non-Arch impies NC}). Taking into account its great
importance in this manuscript, we enunciate here the complete statement.

\begin{corollary}
Let $\left(  X,M,\ast\right)  $ be a non-Archimedean KM-FMS endowed with a
transitive binary relation $\mathcal{S}$ and let $T:X\rightarrow X$ be an
$\mathcal{S}$-nondecreasing fuzzy ample spectrum contraction with respect to
$\left(  M,\mathcal{S},\theta\right)  $. Suppose that $\ast$ is continuous at
the $1$-boundary, $T(X)$ is $(\mathcal{S},d)$-strictly-increasing-precomplete
and there exists a point $x_{0}\in X$ such that $x_{0}\mathcal{S}Tx_{0}$. Also
assume that, at least, one of the following conditions is fulfilled:

\begin{description}
\item[$(a)$] $T$ is $\mathcal{S}$-strictly-increasing-continuous.

\item[$(b)$] $\left(  X,M\right)  $ is $\mathcal{S}$%
-strictly-increasing-regular and condition $(\mathcal{F}_{5})$ holds.

\item[$(c)$] $\left(  X,M\right)  $ is $\mathcal{S}$%
-strictly-increasing-regular and $\theta\left(  t,s\right)  \leq t-s$ for all
$\left(  t,s\right)  \in B\cap(\mathbb{I}\times\mathbb{I})$.
\end{description}

Then the Picard sequence of $T~$based on $x_{0}$ converges to a fixed point of
$T$ (in particular, $T$ has at least a fixed point).

In addition to this, assume that property $(\mathcal{F}_{2}^{\prime})$ is
fulfilled and that for all $x,y\in\operatorname*{Fix}(T)$, there exists $z\in
X$ which is $\mathcal{S}$-comparable, at the same time, to $x$ and to $y$.
Then $T$ has a unique fixed point.
\end{corollary}

We highlight that the previous results also hold in GV-FMS (but we avoid its
corresponding statements).

\subsection{\textbf{Type-2 fuzzy ample spectrum contractions}}

As we commented in Remark \ref{K38 - 32 rem Km to GV}, Definition
\ref{definition KM-space} is as general that such class of fuzzy spaces can
verify that%
\[
M(x,y,t)=0\text{\quad for all }t>0\text{\quad when }x\neq y,
\]
which correspond to infinite distance. In such cases, although a sequence
$\{x_{n}\}$ satisfies%
\[
\theta\left(  M\left(  x_{n+1},x_{n+2},t\right)  ,M\left(  x_{n}%
,x_{n+1},t\right)  \right)  \geq0\quad\text{for all }n\in\mathbb{N}\text{ and
all }t>0,
\]
it is impossible to deduce that $\lim_{n\rightarrow\infty}M\left(
x_{n},x_{n+1},t\right)  =1$ for all $t>0$. Therefore, in order to cover the
fixed point theorems that were demonstrated under this assumption, we must
slightly modify the conditions that a fuzzy ample spectrum contraction
satisfies. In this case, the following new type of contractions may be considered.

\begin{definition}
\label{K38 - 53 def type-2 FASC}Let $\left(  X,M\right)  $ be a fuzzy space,
let $T:X\rightarrow X$ be a self-mapping, let $\mathcal{S}$ be a binary
relation on $X$, let $B\subseteq\mathbb{R}^{2}$ be a subset and let
$\theta:B\rightarrow\mathbb{R}$ be a function. We will say that
$T:X\rightarrow X$ is a \emph{type-2 fuzzy ample spectrum contraction w.r.t.
}$\left(  M,\mathcal{S},\theta\right)  $ if it satisfies properties
$(\mathcal{F}_{1})$, $(\mathcal{F}_{3})$ and the following ones:

\begin{description}
\item[$(\widetilde{\mathcal{F}}_{2})$] If $\{x_{n}\}\subseteq X$ is a Picard
$\mathcal{S}$-strictly-increasing sequence of $T$ such that, for all
$n\in\mathbb{N}$ and all $t>0$,%
\[
M\left(  x_{n},x_{n+1},t\right)  >0\quad\text{and}\quad\theta\left(  M\left(
x_{n+1},x_{n+2},t\right)  ,M\left(  x_{n},x_{n+1},t\right)  \right)  \geq0,
\]
then $\lim_{n\rightarrow\infty}M\left(  x_{n},x_{n+1},t\right)  =1$ for all
$t>0$.

\item[$(\widetilde{\mathcal{F}}_{4})$] If {$x,y\in X$ are such that
$x\mathcal{S}^{\ast}y$ and $Tx\mathcal{S}^{\ast}Ty$ and }$t_{0}>0$ is such
that $M(x,y,t_{0})>0$, then%
\[
{\theta\left(  M(Tx,Ty,t_{0}),M(x,y,t_{0})\right)  \geq0.}%
\]

\end{description}
\end{definition}

In some cases, we will also consider the following properties.

\begin{description}
\item[$(\widetilde{\mathcal{F}}_{2}^{\prime})$] If $\{x_{n}\},\{y_{n}%
\}\subseteq X$ are two $T$-Picard sequences such that, for all $n\in\mathbb{N}
$ and all $t>0$,
\[
x_{n}\mathcal{S}^{\ast}y_{n},\quad{M(x_{n},y_{n},t)>0}\quad\text{and}%
\quad{\theta(M\left(  x_{n+1},y_{n+1},t\right)  ,M\left(  x_{n},y_{n}%
,t\right)  )\geq0}%
\]
then $\lim_{n\rightarrow\infty}{M}\left(  x_{n},y_{n},t\right)  =1$ for all
$t>0$.

\item[$(\widetilde{\mathcal{F}}_{5})$] If $\{\left(  \phi_{n},\psi_{n}\right)
\}$ is a $(T,\mathcal{S}^{\ast},M)$-sequence such that $\phi_{n}(t)>0$,
$\psi_{n}(t)>0$ and $\{\psi_{n}\left(  t\right)  \}\rightarrow1$ for all
$t>0$, and also $\theta(\phi_{n}(t),\psi_{n}(t))\geq0$ for all $n\in
\mathbb{N}$ and all $t>0$, then $\{\phi_{n}\left(  t\right)  \}\rightarrow1$
for all $t>0$.
\end{description}

It is clear that $(\mathcal{F}_{i})\Rightarrow(\widetilde{\mathcal{F}}_{i})$
for all $i\in\{2,4,5\}$ and also $(\mathcal{F}_{2}^{\prime})\Rightarrow
(\widetilde{\mathcal{F}}_{2}^{\prime})$, so each type-1 fuzzy ample spectrum
contraction is a type-2 fuzzy ample spectrum contraction, that is, the notion
of type-2 fuzzy ample spectrum contraction is more general than the notion of
type-1 fuzzy ample spectrum contraction.

\begin{lemma}
\label{K38 - 52 lem equal notion in GV-FMS}In a GV-FMS, the notions of type-1
and type-2 fuzzy ample spectrum contractions coincide.
\end{lemma}

\begin{proof}
It follows from the fact that, in a GV-FMS, $M(x,y,t)>0$ for all $x,y\in X$
and all $t>0$. Hence the respective properties $(\mathcal{F}_{i})$ and
$(\widetilde{\mathcal{F}}_{i})$ are equal.
\end{proof}

However, this generality forces us to assume additional constraints in order
to guarantee existence and uniqueness of fixed points, as we show in the next result.

\begin{theorem}
\label{K38 - 33 th main FAEC type-2}Let $\left(  X,M,\ast\right)  $ be a
KM-FMS endowed with a transitive binary relation $\mathcal{S}$ and let
$T:X\rightarrow X$ be an $\mathcal{S}$-nondecreasing type-2 fuzzy ample
spectrum contraction with respect to $\left(  M,\mathcal{S},\theta\right)  $.
Suppose that $\ast$ is continuous at the $1$-boundary, $T(X)$ is
$(\mathcal{S},d)$-strictly-increasing-precomplete and there exists a point
$x_{0}\in X$ such that $x_{0}\mathcal{S}Tx_{0}$ and $M(x_{0},Tx_{0},t)>0$ for
all $t>0$. Also suppose that the function $\theta$ satisfies:%
\begin{equation}
(t,s)\in B,\quad\theta\left(  t,s\right)  \geq0,\quad s>0\quad\Rightarrow\quad
t>0. \label{K38 - 34 prop theta}%
\end{equation}
Assume that, at least, one of the following conditions is fulfilled:

\begin{description}
\item[$\left(  a\right)  $] $T$ is $\mathcal{S}$-strictly-increasing-continuous.

\item[$\left(  b\right)  $] $\left(  X,M\right)  $ is metrically-$\mathcal{S}%
$-strictly-increasing-regular and condition $(\widetilde{\mathcal{F}}_{5})$ holds.

\item[$\left(  c\right)  $] $\left(  X,M\right)  $ is metrically-$\mathcal{S}%
$-strictly-increasing-regular and $\theta\left(  t,s\right)  \leq t-s$ for all
$\left(  t,s\right)  \in B\cap(\mathbb{I}\times\mathbb{I})$.
\end{description}

If $\left(  X,M\right)  $\ satisfies the property $\mathcal{NC}$, then the
Picard sequence of $T~$based on $x_{0}$ converges to a fixed point of $T$. In
particular, $T$ has at least a fixed point.
\end{theorem}

\begin{proof}
We can repeat many of the arguments showed in the proof of Theorem
\ref{K38 - 31 th main FAEC}, but we must refine them. Let $x_{0}\in X$ be the
point such that $x_{0}\mathcal{S}Tx_{0}$ and $M(x_{0},Tx_{0},t)>0$ for all
$t>0$. Let $\{x_{n}\}$ be the Picard sequence of $T$ starting from $x_{0}$.
Assume that $x_{n}\neq x_{n+1}$ for all $n\in\mathbb{N}_{0}$ and
$x_{n}\mathcal{S}x_{m}$ for all $n,m\in\mathbb{N}_{0}$ such that $n<m$. Since
$T$ is a fuzzy ample spectrum contraction with respect to $\left(
M,\mathcal{S},\theta\right)  $ and $x_{0}\mathcal{S}^{\ast}Tx_{0}=x_{1}$,
$Tx_{0}=x_{1}\mathcal{S}^{\ast}x_{2}=Tx_{1}$ and $M(x_{0},x_{1},t)=M(x_{0}%
,Tx_{0},t)>0$ then%
\[
{\theta\left(  M(x_{1},x_{2},t),M(x_{0},x_{1},t)\right)  ={\theta\left(
M(Tx_{0},Tx_{1},t),M(x_{0},x_{1},t)\right)  }\geq0}%
\]
{for all }$t>0${. Using property (\ref{K38 - 34 prop theta}), it follows that,
for all }$t>0$,%
\[
{\theta\left(  M(x_{1},x_{2},t),M(x_{0},x_{1},t)\right)  \geq0,}\quad
M(x_{0},x_{1},t)>0\quad\Rightarrow\quad M(x_{1},x_{2},t)>0.
\]

By induction, it can be proved that%
\begin{equation}
M(x_{n},x_{n+1},t)>0\quad\text{for all }t>0\text{ and all }n\in\mathbb{N}_{0}
\label{K38 - 35 prop}%
\end{equation}
and%
\[
{\theta\left(  M(x_{n+1},x_{n+2},t),M(x_{n},x_{n+1},t)\right)  \geq0}%
\quad\text{for all }t>0\text{ and all }n\in\mathbb{N}_{0}.
\]
Hence we can apply property $(\widetilde{\mathcal{F}}_{2})$ and we deduce that%
\[
\lim_{n\rightarrow\infty}M\left(  x_{n},x_{n+1},t\right)  \rightarrow
1\quad\text{for all }t>0.
\]
Following the same arguments given in the proof of Theorem
\ref{K38 - 31 th main FAEC}, we can reduce us to the case in which $\{x_{n}\}$
is infinite, in which we know that there is $z\in X$ such that $\{x_{n}\}$
$M$-converges to $z$. It only remain to prove that, under any of conditions
$(a)$, $(b)$ or $(c)$, $z$ is a fixed point of $T$. In case $(a)$, the proof
of Theorem \ref{K38 - 31 th main FAEC} can be repeated.

\begin{description}
\item[$(b)$] Suppose that $\left(  X,M\right)  $ is metrically-$\mathcal{S}%
$-strictly-increasing-regular and condition $(\widetilde{\mathcal{F}}_{5})$
holds. In this case, since $\{x_{n}\}$ is an $\mathcal{S}$-strictly-increasing
sequence such that $\{x_{n}\}\rightarrow z\in X$ and (\ref{K38 - 35 prop})
holds, it follows that%
\[
x_{n}\mathcal{S}z\quad\text{and}\quad M(x_{n},z,t)>0\quad\text{for all }%
n\in\mathbb{N}\text{ and all }t>0.
\]
Taking into account that the sequence $\{x_{n}\}$ is infinite, then there is
$n_{0}\in\mathbb{N}$ such that $x_{n}\neq z$ and $x_{n}\neq Tz$ for all $n\geq
n_{0}-1$. Moreover, as $T$ is $\mathcal{S}$-nondecreasing, then $x_{n+1}%
=Tx_{n}\mathcal{S}Tz$, which means that $x_{n}\mathcal{S}^{\ast}z$ and
$x_{n}\mathcal{S}^{\ast}Tz$ for all $n\geq n_{0}$. Using that $T$ is a type-2
fuzzy ample spectrum contraction with respect to $\left(  M,\mathcal{S}%
,\theta\right)  $ and $M(x_{n},z,t)>0$ for all $n\in\mathbb{N}$ and all $t>0
$, condition $(\widetilde{\mathcal{F}}_{4})$ implies that%
\[
\theta\left(  M(x_{n+1},Tz,t),M(x_{n},z,t)\right)  =\theta\left(
M(Tx_{n},Tz,t),M(x_{n},z,t)\right)  \geq0
\]
for all $n\geq n_{0}$ and all $t>0$. Furthermore, property
(\ref{K38 - 34 prop theta}) ensures that%
\[
\theta\left(  M(x_{n+1},Tz,t),M(x_{n},z,t)\right)  \geq0,\quad M(x_{n}%
,z,t)>0\quad\Rightarrow\quad M(x_{n+1},Tz,t)>0
\]
for all $n\geq n_{0}$ and all $t>0$. Taking into account that $\{M(x_{n}%
,z,t)\}_{n\geq n_{0}}\rightarrow1$ for all $t>0$, assumption $(\widetilde
{\mathcal{F}}_{5})$ applied to the $(T,\mathcal{S}^{\ast},M)$-sequence%
\[
\left\{  \,\left(  \,\phi_{n}=M(Tx_{n},Tz,\cdot),\,\psi_{n}=M(x_{n}%
,z,\cdot)\,\right)  \,\right\}  _{n\geq n_{0}}%
\]
leads to $\{M(x_{n+1},Tz,t)\}_{n\geq n_{0}}\rightarrow1$ for all $t>0$, that
is, $\{x_{n}\}_{n\geq n_{0}}$ $M$-converges to $Tz$. As the $M$-limit in a
KM-FMS is unique, then $Tz=z$, so $z$ is a fixed point of $T$.
\end{description}
\end{proof}

In this context, it is also possible to prove a similar uniqueness result.

\begin{theorem}
\label{K38 - 42 th main FAEC type-2 uniqueness}Under the hypotheses of Theorem
\ref{K38 - 33 th main FAEC type-2}, assume that property $(\mathcal{F}%
_{2}^{\prime})$ is fulfilled and that for all $x,y\in\operatorname*{Fix}(T)$,
there exists $z\in X$ which is $\mathcal{S}$-comparable, at the same time, to
$x$ and to $y$. Then $T$ has a unique fixed point.
\end{theorem}

\begin{proof}
All arguments of the proof of Theorem \ref{K38 - 31 th main FAEC uniqueness}%
\ can be repeated in this context.
\end{proof}

\begin{corollary}
Theorems \ref{K38 - 33 th main FAEC type-2} and
\ref{K38 - 42 th main FAEC type-2 uniqueness} remains true if we replace the
hypothesis that $(X,M)$ satisfies the property $\mathcal{NC}$ by the fact that
$(X,M,\ast)$ is a non-Archimedean KM-FMS.
\end{corollary}

\begin{proof}
It follows from Theorem \ref{K38 - 24 th Non-Arch impies NC}.
\end{proof}

\section{\textbf{Consequences}}

In this section we show some direct consequences of our main results.

\subsection{\textbf{Mihe\textrm{\c{t}}'s fuzzy }$\psi$\textbf{-contractions}}

In \cite{Mi3} Mihe\c{t} introduced a class of contractions in the setting of
KM-fuzzy metric spaces that attracted much attention. It was defined by
considering the following family of auxiliary function. Let $\Psi$ be the
family of all continuous and nondecreasing functions $\psi:\mathbb{I}%
\rightarrow\mathbb{I}$ satisfying $\psi\left(  t\right)  >t$ for all
$t\in\left(  0,1\right)  $. Notice that if $\psi\in\Psi$, then $\psi\left(
0\right)  \geq0$ and $\psi\left(  1\right)  =1$, so $\psi\left(  t\right)
\geq t$ for all $t\in\mathbb{I}$.

\begin{definition}
\label{K38 - 37 def Mihet}\textrm{(Mihe\c{t} \cite{Mi3}, Definition 3.1)}
Given a KM-FMS $\left(  X,M,\ast\right)  $ (it is assumed that $\ast$ is
continuous), a mapping $T:X\rightarrow X$ is a \emph{fuzzy }$\psi
$\emph{-contraction} if there is $\psi\in\Psi$ such that, for all $x,y\in X$
and all $t>0$,%
\begin{equation}
M\left(  x,y,t\right)  >0\quad\Rightarrow\quad M\left(  Tx,Ty,t\right)
\geq\psi\left(  M\left(  x,y,t\right)  \right)  .
\label{K38 - 37 def Mihet, prop}%
\end{equation}

\end{definition}

\begin{theorem}
\label{K38 - 38 th Mihet}\textrm{(Mihe\c{t} \cite{Mi3}, Theorem 3.1)} Let
$(X,M,\ast)$ be an $M$-complete non-Archimedean KM-FMS (it is assumed that
$\ast$ is continuous) and let $T:X\rightarrow X$ be a fuzzy $\psi$-contractive
mapping. If there exists $x\in X$ such that $M(x,Tx,t)>0$ for all $t>0$, then
$T$ has a fixed point.
\end{theorem}

We show that the class of fuzzy ample spectrum contraction properly contains
to the class of Mihe\c{t}'s $\psi$-contractions.

\begin{theorem}
\label{K38 - 39 th Mihet implies type-2}Given a Mihe\c{t}'s fuzzy $\psi
$-contraction $T:X\rightarrow X$ in a KM-FMS $(X,M,\ast)$, let define
$\theta_{\psi}:\mathbb{I}\times\mathbb{I}\rightarrow\mathbb{R}$ by
$\theta_{\psi}\left(  t,s\right)  =t-\psi\left(  s\right)  $ for all
$t,s\in\mathbb{I}$. Then $T$ is a type-2 fuzzy ample spectrum contraction
w.r.t. $(M,\mathcal{S}_{X},\theta_{\psi})$ that also satisfies properties
$(\widetilde{\mathcal{F}}_{2}^{\prime})$, $(\widetilde{\mathcal{F}}_{5})$,
(\ref{K38 - 34 prop theta}) and $\theta_{\psi}\left(  t,s\right)  \leq t-s$
for all $t,s\in\mathbb{I}$.
\end{theorem}

\begin{proof}
Let $T:X\rightarrow X$ be a $\psi$-contraction in a KM-FMS $\left(
X,M,\ast\right)  $. Let consider on $X$ the trivial order $\mathcal{S}_{X}$
defined by $x\mathcal{S}_{X}y$ for all $x,y\in X$. Let $B=\mathbb{I}%
\times\mathbb{I} $ and let define $\theta_{\psi}:\mathbb{I}\times
\mathbb{I}\rightarrow\mathbb{R}$ by $\theta_{\psi}\left(  t,s\right)
=t-\psi\left(  s\right)  $ for all $t,s\in\mathbb{I}$. We claim that $T$ is a
fuzzy ample spectrum contraction w.r.t. $\left(  M,\mathcal{S}_{X}%
,\theta_{\psi}\right)  $ that also satisfies properties $(\widetilde
{\mathcal{F}}_{2}^{\prime})$, $(\widetilde{\mathcal{F}}_{5})$,
(\ref{K38 - 34 prop theta}) and $\theta_{\psi}\left(  t,s\right)  \leq t-s$
for all $t,s\in\mathbb{I}$. We check all conditions.

(\ref{K38 - 34 prop theta}) If $t,s\in\mathbb{I}$ are such that $\theta_{\psi
}\left(  t,s\right)  \geq0$ and $s>0$, then $t-\psi\left(  s\right)  \geq0$,
so $t\geq\psi\left(  s\right)  \geq s>0$. In fact, since $\psi\left(
s\right)  \geq s$ for all $s\in\mathbb{I}$, then%
\[
\theta_{\psi}\left(  t,s\right)  =t-\psi\left(  s\right)  \leq t-s\quad
\text{for all }t,s\in\mathbb{I}.
\]

$(\mathcal{F}_{1})$. It is obvious because $B=\mathbb{I}\times\mathbb{I}$.

$(\widetilde{\mathcal{F}}_{4})$. Let{\ $x,y\in X$ be two points such that
$x\mathcal{S}_{X}^{\ast}y$ and $Tx\mathcal{S}_{X}^{\ast}Ty$ and let }$t>0$ be
such that{\ }$M(x,y,t)>0${. Since }$T$ is a fuzzy $\psi$-contraction, {by
(\ref{K38 - 37 def Mihet, prop}),}%
\[
M\left(  x,y,t\right)  >0\quad\Rightarrow\quad M\left(  Tx,Ty,t\right)
\geq\psi\left(  M\left(  x,y,t\right)  \right)  .
\]
Therefore%
\[
\theta_{\psi}\left(  M\left(  Tx,Ty,t\right)  ,M\left(  x,y,t\right)  \right)
=M\left(  Tx,Ty,t\right)  -\psi\left(  M\left(  x,y,t\right)  \right)  \geq0.
\]

$(\widetilde{\mathcal{F}}_{2}^{\prime})$. Let $x_{1},x_{2}\in X$ be two points
such that, for all $n\in\mathbb{N}$ and all $t>0$,%
\[
T^{n}x_{1}\mathcal{S}^{\ast}T^{n}x_{2},\quad{M(T^{n}x_{1},T^{n}x_{2}%
,t)>0}\quad\text{and}\quad{\theta}_{\psi}{(M(T^{n+1}x_{1},T^{n+1}%
x_{2},t),M(T^{n}x_{1},T^{n}x_{2},t))\geq0.}\quad
\]
Therefore%
\begin{align*}
0  &  \leq{\theta}_{\psi}{(M(T^{n+1}x_{1},T^{n+1}x_{2},t),M(T^{n}x_{1}%
,T^{n}x_{2},t))}\\
&  ={M(T^{n+1}x_{1},T^{n+1}x_{2},t)-\psi(M(T^{n}x_{1},T^{n}x_{2},t)).}%
\end{align*}
Hence ${\psi(M(T^{n}x_{1},T^{n}x_{2},t))\leq M(T^{n+1}x_{1},T^{n+1}x_{2},t)}
$, which means that, for all $n\in\mathbb{N}$ and all $t>0$,%
\begin{equation}
{0<M(T^{n}x_{1},T^{n}x_{2},t)\leq\psi(M(T^{n}x_{1},T^{n}x_{2},t))\leq
M(T^{n+1}x_{1},T^{n+1}x_{2},t)\leq1.} \label{K38 - 40 prop}%
\end{equation}
As a consequence, for each $t>0$, the sequence $\{{M(T^{n}x_{1},T^{n}x_{2}%
,t)}\}_{n\in\mathbb{N}}$ is nondecreasing and bounded above. Hence, it is
convergent. If $L(t)=\lim_{n\rightarrow\infty}{M(T^{n}x_{1},T^{n}x_{2},t)}$,
letting $n\rightarrow\infty$ in (\ref{K38 - 40 prop}) and taking into account
that $\psi$ is continuous, we deduce that%
\[
L(t)\leq\psi\left(  L\left(  t\right)  \right)  \leq L(t).
\]
This means that $\psi\left(  L\left(  t\right)  \right)  =L\left(  t\right)
$, so $L\left(  t\right)  \in\{0,1\}$. However, as $L(t)\geq{M(T^{n}%
x_{1},T^{n}x_{2},t)}>0$, then $L(t)=1$. This proves that $\lim_{n\rightarrow
\infty}{M(T^{n}x_{1},T^{n}x_{2},t)}=1$ for all $t>0$.

$(\widetilde{\mathcal{F}}_{2})$. It follows from $(\widetilde{\mathcal{F}}%
_{2}^{\prime})$ using $x_{2}=Tx_{1}$.

$(\widetilde{\mathcal{F}}_{3})$. Let $\{\left(  \phi_{n},\psi_{n}\right)  \}$
be a $(T,\mathcal{S}^{\ast},M)$-sequence and let $t_{0}>0$ be such that
$\{\phi_{n}\left(  t_{0}\right)  \}$ and $\{\psi_{n}\left(  t_{0}\right)  \}$
converge to the same limit $L\in\mathbb{I}$, which satisfies $L>\phi
_{n}\left(  t_{0}\right)  $ and $\theta_{\psi}(\phi_{n}\left(  t_{0}\right)
,\psi_{n}\left(  t_{0}\right)  )\geq0$ for all $n\in\mathbb{N}$. Therefore,
for all $n\in\mathbb{N}$,%
\[
0\leq\theta_{\psi}(\phi_{n}\left(  t_{0}\right)  ,\psi_{n}\left(
t_{0}\right)  )=\phi_{n}\left(  t_{0}\right)  -\psi\left(  \psi_{n}\left(
t_{0}\right)  \right)  ,
\]
so%
\[
\psi\left(  \psi_{n}\left(  t_{0}\right)  \right)  \leq\phi_{n}\left(
t_{0}\right)  \quad\text{for all }n\in\mathbb{N}.
\]
Letting $n\rightarrow\infty$, as $\psi$ is continuous, we deduce that
$\psi\left(  L\right)  \leq L$, so $L\in\{0,1\}$. However, as $L>\phi
_{n}\left(  t_{0}\right)  \geq0$, then $L=1$.

$(\widetilde{\mathcal{F}}_{5})$. It follows from the fact that $\theta_{\psi
}\left(  t,s\right)  =t-\psi\left(  s\right)  \leq t-s$ for all $t,s\in
\mathbb{I}$.
\end{proof}

The previous result means that each Mihe\c{t}'s fuzzy $\psi$-contraction in a
KM-FMS is a type-2 fuzzy ample spectrum contraction that also satisfies
properties $(\widetilde{\mathcal{F}}_{2}^{\prime})$, $(\widetilde{\mathcal{F}%
}_{5})$, (\ref{K38 - 34 prop theta}) and $\theta_{\psi}\left(  t,s\right)
\leq t-s$ for all $t,s\in\mathbb{I}$. Furthermore, we are going to show that
it is also $\mathcal{S}_{X}$-strictly-increasing-continuous.

\begin{lemma}
\label{K38 - 48 lem Mihet contraction is nondecreasing}Each Mihe\c{t}'s fuzzy
$\psi$-contraction $T:X\rightarrow X$ in a KM-FMS $(X,M,\ast)$ is $T$ is
$\mathcal{S}_{X}$-strictly-increasing-continuous.
\end{lemma}

\begin{proof}
Let $\{x_{n}\}\subseteq X$ be an $\mathcal{S}_{X}$-strictly-increasing
sequence such that $\{x_{n}\}$ $M$-converges to $z\in X$. Let $t_{0}>0$ be
arbitrary. Since $\lim_{n\rightarrow\infty}M(x_{n},z,t_{0})=1$, then there is
$n_{0}\in\mathbb{N}$ such that $M(x_{n},z,t_{0})>0$ for all $n\geq n_{0}$.
Theorem \ref{K38 - 39 th Mihet implies type-2} ensures that $T$ is a type-2
fuzzy ample spectrum contraction w.r.t. $(M,\mathcal{S}_{X},\theta_{\psi})$,
where $\theta_{\psi}:\mathbb{I}\times\mathbb{I}\rightarrow\mathbb{R}$ is given
by $\theta_{\psi}\left(  t,s\right)  =t-\psi\left(  s\right)  $ for all
$t,s\in\mathbb{I}$. Applying property $(\widetilde{\mathcal{F}}_{4})$, we
deduce that, for all $n\geq n_{0}$,%
\[
0\leq\theta_{\psi}(M(Tx_{n},Tz,t_{0}),M(x_{n},z,t_{0}))=M(Tx_{n}%
,Tz,t_{0})-\psi(M(x_{n},z,t_{0})).
\]
As a result, $\psi(M(x_{n},z,t_{0}))\leq M(Tx_{n},Tz,t_{0})\leq1$ for all
$n\geq n_{0}$. Letting $n\rightarrow\infty$, we deduce that%
\[
1=\psi(1)=\lim_{n\rightarrow\infty}\psi(M(x_{n},z,t_{0}))\leq\lim
_{n\rightarrow\infty}M(Tx_{n},Tz,t_{0})\leq1,
\]
so $\lim_{n\rightarrow\infty}M(Tx_{n},Tz,t_{0})=1$. Varying $t_{0}>0$, we
deduce that $\{Tx_{n}\}$ $M$-converges to $Tz$, so $T$ is $\mathcal{S}_{X}$-strictly-increasing-continuous.
\end{proof}

We can conclude that the Mihe\c{t}'s theorem
\ref{K38 - 39 th Mihet implies type-2} is a particular case of our main results.

\begin{corollary}
Theorem \ref{K38 - 38 th Mihet} follows from Theorems
\ref{K38 - 33 th main FAEC type-2} and
\ref{K38 - 42 th main FAEC type-2 uniqueness}.
\end{corollary}

\begin{proof}
Let $(X,M,\ast)$ be an $M$-complete non-Archimedean KM-FMS (it is assumed that
$\ast$ is continuous) and let $T:X\rightarrow X$ be a fuzzy $\psi$-contractive
mapping. Suppose that there exists $x_{0}\in X$ such that $M(x_{0}%
,Tx_{0},t)>0$ for all $t>0$. Given $\psi\in\Psi$, let define $\theta_{\psi
}\left(  t,s\right)  =t-\psi\left(  s\right)  $ for all $t,s\in\mathbb{I}$.
Theorem \ref{K38 - 39 th Mihet implies type-2} guarantees that $T$ type-2
fuzzy ample spectrum contraction w.r.t. $(M,\mathcal{S}_{X},\theta_{\psi})$
that also satisfies properties $(\widetilde{\mathcal{F}}_{2}^{\prime})$,
$(\widetilde{\mathcal{F}}_{5})$, (\ref{K38 - 34 prop theta}) and $\theta
_{\psi}\left(  t,s\right)  \leq t-s$ for all $t,s\in\mathbb{I}$. Notice that
$T(X)$ is $(\mathcal{S}_{X},d)$-strictly-increasing-precomplete because $X$ is
$M$-complete. By Lemma \ref{K38 - 48 lem Mihet contraction is nondecreasing},
$T$ is $\mathcal{S}_{X}$-strictly-increasing-continuous. Item $(a)$ of Theorem
\ref{K38 - 33 th main FAEC type-2} demonstrates that $T$ has, at least, one
fixed point. Furthermore, as $(\widetilde{\mathcal{F}}_{2}^{\prime})$ holds
and any two fixed points of $T$ are $\mathcal{S}_{X}$-comparable, then Theorem
\ref{K38 - 42 th main FAEC type-2 uniqueness} implies that $T$ has a unique
fixed point.
\end{proof}

\begin{example}
Mihe\c{t}'s $\psi$-contractions include Radu's contractions \cite{Radu}
(which, at the same time, generalize Gregori and Sapena's contractions
\cite{GrSa}) satisfying:%
\[
M(Tx,Ty,t)\geq\frac{M(x,y,t)}{M(x,y,t)+k(1-M(x,y,t))}\quad\text{for all
}x,y\in X\text{ and all }t>0,
\]
where $k\in\left(  0,1\right)  $. Therefore, if we consider the function
$\theta:\mathbb{I}\times\mathbb{I}\rightarrow\mathbb{R}$ given, for all
$t,s\in\mathbb{I}$, by%
\[
\theta\left(  t,s\right)  =t-\frac{s}{s+k\left(  1-s\right)  },
\]
then each Radu's contraction is also a fuzzy ample spectrum contraction
associated to $\theta$.
\end{example}

\begin{remark}
In \cite{Mi3} Mihe\c{t} posed the following question: \emph{Does Theorem 3.1
remain true if \textquotedblleft non-Archimedean fuzzy metric
space\textquotedblright\ is replaced by \textquotedblleft fuzzy metric
space\textquotedblright?} Here we cannot give a general answer, but we can say
that if $(X,M)$ satisfies the property $\mathcal{NC}$ and $\ast$ is continuous
at the $1$-boundary, then his Theorem 3.1 remains true.
\end{remark}


\begin{remark}
In order to conclude the current subsection, we wish to highlight one of the main
advantages of the contractivity condition described in $(\mathcal{F}_{4})$,
that is,%
\[
{\theta\left(  M(Tx,Ty,t),M(x,y,t)\right)  \geq0}%
\]
\emph{versus} the Mihe\c{t}'s inequality%
\[
M\left(  Tx,Ty,t\right)  \geq\psi\left(  M\left(  x,y,t\right)  \right)  .
\]
To explain it, let us call $u=M\left(  x,y,t\right)  $ and $v=M\left(
Tx,Ty,t\right)  $. Then the first inequality is equivalent to
\[
{\theta\left(  u,v\right)  \geq0,}%
\]
and the second one can be expressed as%
\[
v\geq\psi(u).
\]
In this sense, Mihe\c{t}'s inequality can be interpreted as a condition in
separate variables, that is, $u$ and $v$ are placed on distinct sides of the
inequality. However, from the inequality ${\theta\left(  u,v\right)  \geq0}$,
in general, it is impossible to deduce a relationship in separate variables
such as $v\geq\psi(u)$. As a consequence, it is often easy to check that a
self-mapping $T:X\rightarrow X$ satisfies a general condition such as
${\theta\left(  u,v\right)  \geq0}$ rather than a more restricted
contractivity condition such as $v\geq\psi(u)$.

To illustrate this advantage, we show how canonical examples of contractions
in the setting of fixed point theory, that is, Banach's contractions in metric
spaces, can be easily seen as fuzzy ample spectrum contractions but, in
general, it is complex to prove that they are also Mihe\c{t}'s fuzzy $\psi$-contractions.
\end{remark}

\begin{example}
Let consider the fuzzy metric $M^{d}$ on $X=\left[  0,\infty\right)  $
associated to the Euclidean metric $d(x,y)=\left\vert x-y\right\vert $ for all
$x\in X$. Given $\lambda\in\left(  0,1\right)  $, let $T:X\rightarrow X$ be
the self-mapping given by%
\[
Tx=\lambda\ln\left(  1+x\right)  \qquad\text{for all }x\in X\text{.}%
\]
Although $T$ is a Banach's contraction associated to the constant $\lambda$,
that is, $d(Tx,Ty)\leq\lambda d(x,y)$ for all $x,y\in X$, in general, it is
not easy to prove that $T$ is a Mihe\c{t}'s $\psi$-contraction because we must
to determine a function $\psi:\mathbb{I}\rightarrow\mathbb{I}$, satisfying
certain properties, and also such that%
\[
\frac{t}{t+\lambda\left\vert \, \ln\dfrac{1+x}{1+y} \, \right\vert }\geq\psi\left(
\frac{t}{t+\left\vert x-y\right\vert }\right)  \qquad\text{for all }x,y\in
X\text{ and all }t>0\text{.}%
\]
To show that $T$ is a fuzzy contraction in $(X,M^{d})$ it could be better to
employ other methodologies, involving terms such that $M^{d}(x,y,\lambda t)$,
rather than Mihe\c{t}'s procedure. However, by handling the function $\theta$
defined as%
\[
\theta\left(  t,s\right)  =\lambda\frac{1-s}{s}-\frac{1-t}{t}\qquad\text{for
all }t,s\in\left(  0,1\right]  ,
\]
it can directly checked that $T$ is a fuzzy ample spectrum contraction
because, for all $x,y\in X$ and all $t>0$:%
\begin{align*}
\theta(M^{d}\left(  Tx,Ty,t\right)  ,M^{d}\left(  x,y,t\right)  )  =
\lambda \, \frac{~1-\dfrac{t}{t+d(x,y)}~}{\dfrac{t}{t+d(x,y)}}-\frac{~1-\dfrac
{t}{t+d(Tx,Ty)}~}{\dfrac{t}{t+d(Tx,Ty)}} =\frac{\lambda d(x,y)-d(Tx,Ty)}{t}.
\end{align*}
\end{example}

\subsection{\textbf{Altun and Mihe\textrm{\c{t}}'s fuzzy contractions}}

In \cite{AlMi} Altun and Mihe\c{t} introduce the following kind of fuzzy
contractions in the setting of ordered KM-FMS. Recall that a \emph{partial
order on a set }$X$ is a binary relation $\mathcal{S}$ on $X$ which is
reflexive, antisymmetric and transitive.

\begin{theorem}
\textrm{(Altun and Mihe\c{t} \cite{AlMi}, Theorem 2.4)} Let $(X,M,\ast)$ be a
complete non-Archimedean KM-FMS (it is assumed that $\ast$ is continuous) and
let $\preceq$ be a partial order on $X$. Let $\psi:\mathbb{I}\rightarrow
\mathbb{I}$ be a continuous mapping such that $\psi(t)>t$ for all $t\in\left(
0,1\right)  $. Also let $T:X\rightarrow X$ be a nondecreasing mapping w.r.t.
$\preceq$ with the property%
\[
M(Tx,Ty,t)\geq\psi(M(x,y,t))\quad\text{for all }t>0\text{ and all }x,y\in
X\text{ such that }x\preceq y.
\]
Suppose that either:

\begin{description}
\item[$(a)$] $T$ is continuous,\quad or

\item[$(b)$] $x_{n}\preccurlyeq x$ for all $n\in\mathbb{N}$ whenever
$\{x_{n}\}\subseteq X$ is a nondecreasing sequence with $\{x_{n}\}\rightarrow
x\in X$
\end{description}

hold. If there exists $x_{0}\in X$ such that%
\[
x_{0}\preccurlyeq Tx_{0}\quad\text{and}\quad\psi(M(x_{0},Tx_{0},t))>0\quad
\text{for each }t>0,
\]
then $T$ has a fixed point.
\end{theorem}

Obviously, the previous theorem can be interpreted as a version of Theorem
\ref{K38 - 37 def Mihet} in which a partial order (which is a transitive
binary relation) is employed. Notice that here the function $\psi$ is not
necessarily nondecreasing, but such property have not been used in the
arguments of the proofs of the previous subsection.

\subsection{\textbf{Fuzzy ample spectrum contractions by using admissible mappings}}

{Following \cite{SaVeVe}, }when a fuzzy space $(X,M)$ is endowed with a
function $\alpha:X\times X\rightarrow\left[  0,\infty\right)  $, the following
notions can be introduced:

\begin{itemize}
\item the function $\alpha$ is \emph{transitive} if $\alpha(x,z)\geq1$ for
each point $x,z\in X$ for which there is $y\in X$ such that $\alpha(x,y)\geq1$
and $\alpha(y,z)\geq1$;

\item a mapping $T:X\rightarrow X$ is $\alpha$\emph{-admissible} whenever
$\alpha\left(  x,y\right)  \geq1$ implies $\alpha\left(  Tx,Ty\right)  \geq1$,
where $x,y\in X$;

\item a sequence $\{x_{n}\}\subseteq X$ is $\alpha$\emph{-nondecreasing} if
$\alpha(x_{n},x_{n+1})\geq1$ for all $n\in\mathbb{N}$;

\item a sequence $\{x_{n}\}\subseteq X$ is $\alpha$\emph{-strictly-increasing}
if $\alpha(x_{n},x_{n+1})\geq1$ and $x_{n}\neq x_{n+1}$ for all $n\in
\mathbb{N}$;

\item a mapping $T:X\rightarrow X$ is $\alpha$\emph{-nondecreasing-continuous}
if $\{Tx_{n}\}\rightarrow Tz$ for all $\alpha$-nondecreasing sequence
$\{x_{n}\}\subseteq X$ such that $\{x_{n}\}\rightarrow z\in X$;

\item a mapping $T:X\rightarrow X$ is $\alpha$%
\emph{-strictly-increasing-continuous} if $\{Tx_{n}\}\rightarrow Tz$ for all
$\alpha$-strictly-increasing sequence $\{x_{n}\}\subseteq X$ such that
$\{x_{n}\}\rightarrow z\in X$;

\item a subset $Y\subseteq X$ is $(\alpha,M)$%
\emph{-strictly-increasing-complete} if every $\alpha$-strictly-increasing and
$M$-Cauchy sequence $\{y_{n}\}\subseteq Y$ is $M$-convergent to a point of $Y$;

\item a subset $Y\subseteq X$ is $(\alpha,M)$%
\emph{-strictly-increasing-precomplete} if there exists a set $Z$ such that
$Y\subseteq Z\subseteq X$ and $Z$ is $(\alpha,M)$-strictly-increasing-complete;

\item $\left(  X,M\right)  $ is $\alpha$\emph{-strictly-increasing-regular}
if, for all $\alpha$-strictly-increasing sequence $\{x_{n}\}\subseteq X$ such
that $\{x_{n}\}\rightarrow z\in X$, it follows that $\alpha(x_{n},z)\geq1$ for
all $n\in\mathbb{N}$;
\end{itemize}

In this setting, it is possible to introduce the notion of \emph{fuzzy ample
spectrun contraction w.r.t. }$(M,\alpha,\theta)$ by replacing any condition of
the type $x\mathcal{S}y$ by $\alpha(x,y)\geq1$ in Definition
\ref{K38 - 54 def FASC} (or Definition \ref{K38 - 53 def type-2 FASC}). In
this general framework, it is possible to obtain some consequences as the
following ones.

\begin{corollary}
Let $\left(  X,M,\ast\right)  $ be a KM-FMS endowed with a transitive function
$\alpha:X\times X\rightarrow\left[  0,\infty\right)  $ and let $T:X\rightarrow
X$ be an $\alpha$-nondecreasing fuzzy ample spectrum contraction with respect
to $\left(  M,\alpha,\theta\right)  $. Suppose that $\ast$ is continuous at
the $1$-boundary, $T(X)$ is $(\alpha,d)$-strictly-increasing-precomplete and
there exists a point $x_{0}\in X$ such that $\alpha(x_{0},Tx_{0})\geq1$. Also
assume that, at least, one of the following conditions is fulfilled:

\begin{description}
\item[$(a)$] $T$ is $\alpha$-strictly-increasing-continuous.

\item[$(b)$] $\left(  X,M\right)  $ is $\alpha$-strictly-increasing-regular
and condition $(\mathcal{F}_{5})$ holds.

\item[$(c)$] $\left(  X,M\right)  $ is $\alpha$-strictly-increasing-regular
and $\theta\left(  t,s\right)  \leq t-s$ for all $\left(  t,s\right)  \in
B\cap(\mathbb{I}\times\mathbb{I})$.
\end{description}

If $\left(  X,M\right)  $\ satisfies the property $\mathcal{NC}$, then the
Picard sequence of $T~$based on $x_{0}$ converges to a fixed point of $T$. In
particular, $T$ has at least a fixed point.
\end{corollary}

\begin{proof}
It follows from Theorem \ref{K38 - 31 th main FAEC} by considering on $X$ the
binary relation $\mathcal{S}_{\alpha}$ given, for each $x,y\in X$, by
$x\mathcal{S}_{\alpha}y$ if $\alpha(x,y)\geq1$.
\end{proof}

\vspace*{-3mm}

\section{\textbf{Conclusions and prospect work}}

In this paper we have introduced the notion of \emph{fuzzy ample spectrum
contraction} in the setting of fuzzy metric spaces in the sense of Kramosil
and Mich\'{a}lek into two distinct ways. After that, we have demonstrated some
very general theorems in order to ensure the existence and uniqueness of fixed
points for such families of fuzzy contractions. We have also illustrated that
these novel classes of fuzzy contractions extend and generalize some well
known groups of previous fuzzy contractions.

In order to attract the attention of many researchers in the field of fixed
point theory, in the title we announced that out results was going to be
developed in the setting of non-Archimedean fuzzy metric spaces. However, we
have presented a new property (that we have called $\mathcal{NC}$) in order to
consider more general families of fuzzy metric spaces. By working with
property $\mathcal{NC}$ we have given a positive partial answer to an open
problem posed by Mihe\c{t} in \cite{Mi3}.

In this line of research there is much prospect work to carry out. For
instance, we have shown that these novel fuzzy contraction generalize the
notion of \emph{ample spectrum contraction} in the setting of metric spaces
under an additional assumption. Then, the following doubts naturally appears:

\emph{Open question 1:} is any ample spectrum contraction in a metric space
$(X,d)$ a fuzzy ample spectrum contraction in the fuzzy metric space
$(X,M^{d})$?

\emph{Open question 2:} in order to cover other kinds of fuzzy contractions,
is it possible to extend the notion of ample spectrum contraction to the fuzzy
setting by introducing a function in the argument $t$ of $M(x,y,t)$?

\vspace*{3mm}

\section*{Acknowledgments}

A.F. Rold\'{a}n L\'{o}pez de Hierro is grateful to Project TIN2017-89517-P of
Ministerio de Econom\'{\i}a, Industria y Competitividad and also to Junta de
Andaluc\'{\i}a by project FQM-365 of the Andalusian CICYE. This article was
funded by the Deanship of Scientific Research (DSR), King Abdulaziz
University, Jeddah. N. Shahzad acknowledges with thanks DSR for financial support.

\end{document}